\documentclass[]{interact}

\usepackage{epstopdf}
\usepackage{subfigure}
\usepackage{lscape}
\usepackage{multicol} 
\usepackage{multirow} 

\usepackage{natbib}
\bibpunct[, ]{(}{)}{;}{a}{}{,}
\renewcommand\bibfont{\fontsize{10}{12}\selectfont}

\theoremstyle{plain}
\newtheorem{theorem}{Theorem}[section]

\newtheorem{proposition}[theorem]{Proposition}

\theoremstyle{definition}

\theoremstyle{remark}

\begin{document}

\articletype{ARTICLE TEMPLATE}

\title{A fresh look at symmetric traffic assignment and algorithm convergence}

\author{
\name{Priyadarshan N. Patil\textsuperscript{a}\thanks{CONTACT Priyadarshan N. Patil. Email: Priyadarshan@utexas.edu}}
\affil{\textsuperscript{a}Graduate Program in Operations Research and Industrial Engineering, The University of Texas at Austin, 204 E. Dean Keeton St., C2200, Austin, TX, USA.}
}

\maketitle

\begin{abstract}
Extensions of the static traffic assignment problem with link interactions were studied extensively in the past. Much of the network modeling community has since shifted to dynamic traffic assignment incorporating these interactions. We believe there are several reasons to re-examine static assignment with link interactions. First, if link interactions can be captured in a symmetric, monotone manner, equilibrium always exists and is unique, and provably-correct algorithms exist. We show that several of the most efficient algorithms for the separable traffic assignment problem can be readily applied with symmetric interactions. We discuss how the (asymmetric) Daganzo merge model can be approximated by symmetric linear cost functions. Second, we present computational evidence suggesting that convergence to equilibrium is \emph{faster} when symmetric, monotone link interactions are present. This is true even when interactions are asymmetric, despite the lack of a provable convergence result. Lastly, we present convergence behavior analysis for commonly used network and link metrics. For these reasons, we think static assignment with link interactions deserves additional attention in research and practice.

\end{abstract}


\section{Introduction}
\label{sec:introduction}

Traffic assignment is a cornerstone of the urban planning and forecasting process, relating travel demand to a forecasted network loading, providing metrics such as link flows and travel times. Many traffic assignment models exist. The most common in practice is the static traffic assignment problem (TAP) first formulated by \cite{Beckmann}. Dynamic traffic assignment has been researched for several decades, and is increasingly used by practitioners. Both static and dynamic traffic assignment are described at length in~\cite{blubook_vol1_v085}. Equilibrium assignment is also now a feature of some microsimulation software. Broadly speaking, there is a tension between the level of ``realism'' of a model (used here to refer to the level of detail in capturing traffic physics) and how ``nice'' the model is (in terms of mathematical analysis, computational efficiency, provably correct algorithms, and properties such as equilibrium existence, uniqueness, and stability). Dynamic link and node models are much more descriptive of traffic flow than the volume-delay functions typical of static assignment, a powerful argument in their favor. At the same time, dynamic equilibrium need not exist~\citep[Section 11.3.1]{blubook_vol1_v085}, or several may exist~\citep{nie10eqa}; and precisely the same features that make dynamic models more realistic may also make them more sensitive to any errors in input data, complicating calibration and possibly introducing more error than is saved by improving the traffic flow model~\citep{boyles2019}. These are not merely theoretical concerns --- if dynamic equilibrium may not exist at all, or if multiple equilibria exist, it is not at all clear how to use the results of such a model for alternatives analysis; and if a model is highly sensitive to parameters which are hard to estimate (such as a time-dependent OD matrix decades in the future) its practical utility is heavily limited. More extended versions of these issues are discussed in~\cite{bar2010}, \cite{boyles2019}, and \cite{Patil2020}; in short, some applications demand the realism of a dynamic model, while in other applications the advantages of static models outweigh their drawbacks, and careful researchers and practitioners select their tool based on the problem at hand.

This ongoing conversation serves as background for our investigation.
Some researchers have attempted to find unifying frameworks for static and dynamic assignment; see \cite{bliemer2012quasi} and \cite{bliemer2019static} for examples.
An alternative approach is to improve the traffic model in static assignment.
In the 1980s, there was an active line of research in the static assignment problem with \emph{interactions} among links, rather than using separable link performance functions.
(This background material is described at greater length in Section~\ref{sec:background}.)
While research in this area has continued to this day, much of the community's attention shifted to dynamic assignment in the 1990s, with the advent of reasonable link and node models, such as the cell transmission model~\citep{daganzo94,daganzo95}.
In recent decades, our knowledge of good dynamic link and node models has advanced further still~\citep{yperman_diss, tampere11}.
Research has also advanced considerably in how the separable TAP is solved, with the discovery of path- and bush-based algorithms~\citep{jaya,bar2002,Dial,bar2010,xie,chen2020parallel}.

We believe there are several reasons to take a fresh look at the static traffic assignment problem with symmetric link interactions, considering it at least as an alternative to static assignment, if not dynamic.
As is known, with monotone cost functions, S-TAP retains most of the favorable properties of TAP, including formulating the equilibrium problem as a convex program, and the resulting features of solution existence, uniqueness, and algorithmic tractability.
In this study, we accomplish the following:
\begin{enumerate}
\item As an example of how existing node models can be approximated by symmetric, monotone link performance functions, we develop such a representation of the Jin-Zhang merge model~\citep{jin2003distribution}.  (Section~\ref{sec:merge})
\item We discuss solution algorithms for S-TAP, including classic convex combinations methods, but focusing mainly on more recent algorithms based on shifting flow between a pair of alternative segments. 
We show that the flow shift formula for S-TAP takes a familiar and simple form, and therefore existing algorithms for TAP can be easily adapted for S-TAP.  (Section~\ref{sec:convergence}).
\item We implement these algorithms on standard test networks, and show that in most cases, S-TAP actually converges \emph{faster} than the separable TAP.
Therefore, S-TAP can be considered as a serious alternative to TAP in planning practice.
We also report preliminary results showing good performance of gradient projection for asymmetric link interactions, even though there is no convex objective function (and therefore no guarantees of convergence).
\end{enumerate}
Section~\ref{sec:algconv_conclusions} summarizes our findings, and provides specific suggestions for further investigations on the utility of S-TAP.

\section{Model Formulation}
\label{sec:formulation}
Consider a network with node and link sets $N$ and $A$, respectively.
Let $Z$ denote the subset of nodes where trips may start and end, and for each $(r,s) \in Z^2$, let $d_{rs}$ denote the fixed travel demand from node $r$ to node $s$, and $\Pi_{rs}$ the set of paths connecting $r$ to $s$.
Let $h_\pi$ denote the number of travelers using path $\pi$.
For each link $a \in A$, let $x_a$ and $t_a$ denote its flow and travel time.
If each $t_a$ is a nonnegative, continuous, and increasing function of $x_a$ alone, then the flows $\mathbf{x}$ solving the convex program\newline
\begin{equation}
\label{eqn:integral_defn}
\qquad    \min_{\mathbf{x},\mathbf{h}}  \sum_{a} \int_0^{x_a} t_a (x)~dx
\end{equation}
\text{subject to:}
\begin{align*}
    x_{a} &= \sum_{\pi \in \Pi : a \in \pi} h_\pi & \forall a \in A\\
    \sum_{\pi \in \Pi_{rs}} h_\pi &= d_{rs} & \forall(r,s) \in Z^2\\
    h_\pi &\geq 0 & \forall \pi \in \Pi
\end{align*}

are an equilibrium, with all used paths having equal and minimal travel time~\citep{Beckmann}.
We will use $X$ to denote the set of link flows corresponding to feasible solutions of TAP; this set is compact and convex.
We will also refer to this optimization problem as TAP.
As the objective is strictly convex, equilibrium exists and is unique in link flows.
There are many efficient algorithms for solving TAP, including convex combinations algorithms~\citep{powell1982convergence,frank,chen2002,mitradjieva2013stiff}, path-based algorithms~\citep{larsson1992simplicial,jaya,lee2003,florian,kumar2010,kumar2014,xie2018}, and bush-based algorithms~\citep{bar2002,Dial,gentile,bar2010,nie2012note,xie2}.
Many (but not all) of the latter two algorithm types shift flow between pairs of path segments as a key operation, aiming to equilibrate them using a line search.

Now instead assume that the link performance functions depend on multiple links' flow. For generality, we write $t_a(\mathbf{x})$ to express dependence on (potentially) every other link flow, although in practice each link's travel time depends only on a few other links.
If these functions are differentiable, and the Jacobian of $\mathbf{t}$ with respect to $\mathbf{x}$ is symmetric, equilibrium flows correspond to stationary points of the function 
\begin{equation}
\label{eqn:line_integral_defn}
\qquad    \min_{\mathbf{x},\mathbf{h}}  \oint_0^{\vec{x}} t_i (\vec{s})  d\vec{s}  
\end{equation}

relative to the same constraints as TAP, as can be seen by writing the optimality conditions.
We will refer to this optimization problem as S-TAP.
The symmetry requirement on the Jacobian is critical.
Otherwise, the line integral is path-dependent and the function as stated is not well-defined, and the model would need to be formulated and solved using variational inequality (VI) methods.~\cite{fach2003}
We will refer to the asymmetric case as A-TAP.

If this Jacobian is additionally positive definite, then we say that the link performance functions are \emph{monotone} (note that this is stronger than requiring each $t_a$ to be monotone in each link flow separately).
In such a case, the objective~\eqref{eqn:line_integral_defn} is a strictly convex function (its Hessian is the Jacobian just described), and equilibrium again exists and is unique.
If the link performance functions are not monotone, the function \eqref{eqn:line_integral_defn} is not convex, and multiple stationary points (and therefore equilibria) may exist, with the minima corresponding to stable equilibria.  
For example, see Figure~\ref{fig:S-TAP_multiple_optima}.
There are three equilibria corresponding to the flow vectors (0,10), (5,5), and (10,0); the first and third of these are stable, and the second is unstable.
For derivations of the above results and more discussion, see the books by~\cite{patriksson2015traffic} and~\cite{blubook_vol1_v085}.

\section{Background}
\label{sec:background}

Historically, \cite{prager1954problems} first mentioned the need to model traffic interactions on a two-way street.
\cite{dafermos1971extended} and \cite{dafermos1972traffic} were the first to formulate S-TAP and A-TAP, showing equivalence with the multi-class TAP, and presenting an iterative flow update algorithm to obtain user equilibrium and system optimal flows. 
\cite{sender1970equilibre} used the fixed point theorem to show the existence of an equilibrium solution in the asymmetrically formulated multi-class TAP. 
\cite{smith1979existence} was the first to formulate A-TAP as a VI, presenting a set of conditions for existence (continuity of cost vectors) and uniqueness (strict monotonicity of cost vectors for all supply feasible vectors) of the equilibrium solution.
\cite{dafermos1980traffic} showed that these uniqueness conditions were equivalent to the Jacobian being positive definite, and proved the existence of an equilibrium solution.
As checking positive definiteness can be cumbersome to test in practice, \cite{heydecker1983some} proposed an easier test, based on diagonal dominance.
This is a weaker condition than positive definiteness, thus showing that positive definiteness of the Jacobian is a sufficient but not necessary condition.
Heydecker also discussed the existence of multiple (and unstable) equilibrium solutions when said conditions are violated.
The equivalent requirement of positive definite Jacobian for the multi-class TAP solution existence was also shown by \cite{braess1979existence}.

The VI formulation allowed for a number of different solution approaches to be proposed for A-TAP, including the non-linear Jacobi method (or diagonalization method), projection methods, and column generation methods.
\cite{dafermos1982relaxation} proposed a relaxation method and proved its convergence. 
\cite{fisk1982solution} analyzed this approach, the projection method, and three other solution methods, concluding that the non-linear Jacobi method was most efficient.  
\cite{nagurney1984comparative,nagurney1986computational} extended this comparison of the relaxation and projection method for multi-modal problem with varying travel costs and A-TAP, concluding that no one method was uniformly superior.
\cite{smith1983existence} modeled junction interactions using A-TAP, while proposing a new objective function measuring deviation of traffic distribution from equilibrium.
Smith also provided a descent direction to monotonically reduce the new objective function and a path enumeration based algorithm for A-TAP. 

\cite{florian1982convergence} provided a sufficient condition for diagonalization algorithm convergence.
\cite{nguyen1984efficient} proposed an iterative approach for A-TAP post-optimizing a linear sub-problem at each iteration.
\cite{lawphongpanich1984simplical} used simplicial decomposition, constructing the solution as a convex combination of all extreme points of the flow vector space.
They compared their approach to the \cite{nguyen1984efficient} approach on the networks proposed in ~\cite{nguyen1984efficient} and ~\cite{fisk1982solution}, concluding that the simplicial approach is competitive with the iterative approach while storing a small number of extreme flow patterns.
In the same vein, \cite{gabriel1997traffic} and \cite{bernstein1997solving} solved A-TAP with non-additive route costs using column generation for nonlinear complementarity problem (NCP) gap minimization, equivalent to the VI formulation.
The NCP formulation is needed when route-level interactions are modeled, such as non-additive route cost interactions.
They also showed the existence and uniqueness of the solution if the function is separable.
\cite{lo2000traffic} also apply column generation to the reformulated NCP for route-specific tolls.

\cite{mahmassani1988some} tested the diagonalization approach on three networks including the Texas highway network, and observed convergence of the algorithm despite violation of sufficient conditions presented in~\cite{dafermos1982relaxation}.
Similar observations were made by \cite{friesz1984alternative} and \cite{friesz1985transportation}.
Mahmassani and Mouskos also compared their implementations with Sheffi's streamlined implementation~\cite{sheffi1985urban} with one FW iteration per subproblem and concluded that there is no single best choice for the number of FW iterations.
\cite{meneguzzer1995equilibrium} provided an overview of the advances in the field of diagonalization for A-TAP and convergence for explicit modeling of intersections.
\cite{marcotte1988adaptation} applied the modified Newton method to A-TAP, comparing its performance with cutting plane methods and diagonalization, observing their method performing better than other methods for complex asymmetric interactions. 
\cite{dupuis1986convergence} assessed the convergence conditions for A-TAP solutions using projection and diagonalization methods.

\cite{hearn1984convex} drew a connection between convex programming formulation of A-TAP and the VI formulation.
\cite{marcotte2004new} relaxed the monotonicity condition for interactions, proposing weaker convergence conditions for A-TAP with multiple modes.
\cite{wong2001path} modeled A-TAP with simulation approach for the intersection delay.
\cite{panicucci2007path} formulated the VI in terms of path flows and propose a column generation scheme based on Khobotov's method.
\cite{yook2016acceleration} discussed ways to accelerate the convergence of double-projection method proposed by Panicucci et al. using the decomposable path flow VI structure.
\cite{grange2009equivalent} presented a method to equalize line integral paths for affine cost A-TAP.
\cite{chen2011modeling} modeled A-TAP interactions using side constraints.
\cite{sancho2015applying} evaluated the performance of five projection methods for A-TAP and observe that a variation of Khobotov's method proposed by \cite{he2012proximal} shows the fastest convergence.
\cite{patriksson2015traffic} and \cite{yook2014models} provide good overviews of alternative solution approaches for A-TAP. 

As suggested by the above review, there was significant focus on A-TAP research in the 1980s and 1990s, but comparatively less since then.
The broad explanation for this is the emergence of DTA as a serious modeling tool~\citep{peeta01,chiu2011dynamic}, with better grounding in traffic science than link performance functions could ever hope to have, whether separable or with interactions.  
Nevertheless, as suggested in the introduction, there are other advantages TAP and its variants have, and both tools are used today for different applications.
Our aim in this paper is to suggest that further work in S-TAP and A-TAP can be valuable, incorporating insights from recent DTA research, adapting recent algorithms for the basic TAP, and showing that S-TAP may be a reasonable alternative in places where static assignment is still useful.

\section{A symmetric, monotone merge model}
\label{sec:merge}

While much research has been done on theoretical properties of traffic assignment with interactions, we are not aware of specific guidance on exactly \emph{how} interactions should be chosen to represent real-world scenarios.
To further motivate our investigation, we will give an illustration of how symmetric, monotone link performance functions can be chosen to approximate a node model used in dynamic traffic assignment.
The method we describe here is surely not the only way to do this, but a full study of how to develop such approximations (and to assess their quality on full-sized, realistic networks) is beyond the scope of this study.
Our aim in this section is simply to demonstrate that S-TAP is a plausible model for certain applications.

We take as our starting point a simple network loading model, a network of point queues.
In this model, the time required to travel each link is a fixed free-flow time ($t^0$), plus time spent waiting in a queue at the downstream end.
Let $Q(\tau)$ denote the length of the queue at time $\tau$; these queues are ``point'' queues in that they may grow arbitrarily long.  
Each link also has a uniform saturation flow $u$, perhaps proportional to the number of lanes.
We will denote by $x(\tau)$ the \emph{inflow} rate at the link's upstream end, and $y(\tau)$ the \emph{outflow} rate at the link's downstream end.
If $x$ and $y$ are time-invariant, we must have $x \geq y$, and the queue length at a given point in time will be $Q(\tau) = \tau(x - y)$.
The travel time experienced for a vehicle entering the link at time $t$ will be $t^0 + \tau( x/y - 1)$.
If we assume that $x$ is constant over an interval of length $T$ (and zero otherwise), the average delay experienced by a vehicle on the link will be $t^0 + \frac{T}{2} ( x/y - 1)$.
If we choose units so that $T = 2$, the formula simplifies to $t^0 + (x/y - 1)$, which we will adopt for the remainder of the section.

Consider now a merge node with two upstream links (indexed 1 and 2) and one downstream link (indexed by 3).
A ``merge model'' takes as input the demands from each upstream link and the available supply from the downstream link, and produces the flow rates out of each upstream link into the downstream link.
Several merge models have been proposed in the literature~\citep{daganzo95,lebacque96,jin2003distribution}.

We now show that the Jin-Zhang merge model~\citep{jin2003distribution} leads to a symmetric, monotone delay function, if the demands are interpreted as the inflows $x$.\footnote{This is an approximation to the actual dynamic model, since if a queue forms, the demand will increase to the saturation flow $u$.
Indeed, one critique of the Jin-Zhang model is that it is unstable with respect to this transition, a violation of the ``invariance principle,'' cf.~\cite{lebacque05}.
}
With this interpretation, the model takes the following form.
If $x_1 + x_2 \leq u_3$, there are no queues, and hence no delays: $y_1 = x_1$, $y_2 = x_2$, so $t_1 = t^0_1$ and $t_2 =t^0_2$.
Otherwise, the Jin-Zhang model allocates flows proportionally to demands: $y_1 = u_3 x_1 / (x_1 + x_2)$, and hence $t_1 = t^0_1 + [(x_1 + x_2) / u_3 - 1]$.
Likewise, $t_2 = t^0_2 + [(x_1 + x_2) / u_3 - 1]$.
In either case, the Jacobian of $\mathbf{t}$ with respect to $\mathbf{x}$ is symmetric and positive semidefinite.

Alternative merge models, such as that of \cite{daganzo95}, do not directly lead to symmetric and monotone performance functions.
However, it may be possible to create reasonable approximations to them that satisfy these conditions (indeed the Jin-Zhang model may be seen as such an approximation), at least in a region of demand and supply values the modeler believes to be likely at a particular junction. 
For instance, the Daganzo merge violates symmetry only in the exceptional case when a queue forms on one upstream link, not both.
We believe it plausible that similar procedures or approximations can apply to other types of interactions between links, based either on node models from dynamic traffic assignment, formulas from the \cite{manual2010hcm2010} or similar literature, or regression from simulation, but we leave such investigation to future study.

\section{Algorithms and convergence}
\label{sec:convergence}

This section discusses solution algorithms for S-TAP with monotone cost functions.
To do so, it will be convenient to choose a specific integration path for the line integral in \eqref{eqn:line_integral_defn}.
If we choose the path
$(0, 0, 0, \ldots,0) \rightarrow
(x_1, 0, 0, \ldots, 0) \rightarrow
(x_1, x_2, 0, \ldots, 0) \rightarrow
\cdots \rightarrow
(x_1, x_2, \ldots, x_n)
\,,
$
the line integral decomposes into a sum of ordinary integrals as
\begin{equation}
\label{eqn:line_integral_split}
F(\mathbf{x})
=
\sum_{i=1}^n \int_{(x_1, \ldots, x_{a-1}, 0, 0, \ldots, 0)}^{(x_1, \ldots, x_{a-1}, x_a, 0, \ldots, 0)} t_i(x)~dx
\,.
\end{equation} 
Each link performance function is increasing for each independent flow variable, so each integral is convex, as is the sum.

With this representation, the gradient and Hessian of $F$ take simple forms.
Using the Leibniz rule, the derivative of $F$ with respect to any link flow $x_a$ is
\[
\frac{\partial F}{\partial x_a} =
t_a(x_1, \ldots, x_a, 0, \ldots, 0) + \sum_{i > a} \int_{(x_1, \ldots, x_{i-1}, 0, \ldots, 0)}^{(x_1, \ldots, x_{i-1}, x_{i}, 0, \ldots, 0)} \frac{\partial t_i}{\partial x_a}(x)~dx
\,.
\]
But using the symmetry condition, $\partial t_i/\partial x_a = \partial t_a/\partial x_i$, and so we ultimately have $\partial F / \partial x_a = t_a$, from the fundamental theorem of calculus.
That is, the gradient of $F$ is simply the vector of link travel times, as it is in TAP (and indeed, as it should be for the optimality conditions to express equilibrium).
The Hessian of $F$ is then just the Jacobian of the link performance functions, that is, $(HF)_{ab} = \partial t_a / \partial x_b$.

\subsection{Convex combinations algorithms}

Convex combinations algorithms operate on the link flow vector $\mathbf{x}$, iteratively combining a current feasible solution with a ``target'' solution $\mathbf{x^*}$, with an update rule of the form $\mathbf{x} \leftarrow (1 - \lambda) \mathbf{x} +  \lambda \mathbf{x^*}$.
Typically $\mathbf{x^*}$ is an ``all-or-nothing'' solution obtained by placing all demand on shortest paths when the link costs are $\mathbf{t(x)}$.
They are relatively na\"{i}ve, but amenable to parallelization, and they typically make excellent progress in their initial iterations before the rate of convergence slows sharply (ultimately, to a sublinear rate).

The simplest possible convex combinations algorithm is the method of successive averages, where the step sizes $\lambda$ are chosen \emph{a priori} in a divergent sequence (but with $\sum \lambda^2$ finite); a common example is $ \{ 1/2, 1/3, 1/4, \ldots \}$.
Convergence of this method for S-TAP can be shown using the following result:

\begin{proposition}
\emph{\citep{powell1982convergence}}
\label{prp:powellsheffi}
Consider a twice-continuously differentiable convex function $F(\mathbf{x})$, and a sequence  $\{ \mathbf{x_0}, \mathbf{x_1}, \ldots \}$, where $\mathbf{x_0} \in X$ and $\mathbf{x_i} = (1 - \lambda_i) \mathbf{x_{i-1}} + \lambda_i \mathbf{x^*_i}$ for $i \geq 1$, with $\mathbf{x^*_i} \in X$, $\lambda_i \in [0, 1]$, $\sum \lambda_i = \infty$, but $\sum \lambda_i^2 < \infty$.
This sequence converges to a minimizer $\bar{\mathbf{x}}$ of $F$ if the following conditions hold:
\begin{enumerate}
\item The inner product $(\nabla F (\mathbf{x_i}))^T (\mathbf{x^*_i} - \mathbf{x_i})$ is negative  whenever $F(\mathbf{x_i}) > F(\bar{\mathbf{x}})$.
\item The values of $(\mathbf{x^*_i} - \mathbf{x_i})^T HF(\mathbf{x_i} + \psi \lambda_i (\mathbf{x^*_i} - \mathbf{x_i}) ) (\mathbf{x^*_i} - \mathbf{x_i})$ are bounded over all $i$ and $\psi \in [0, 1]$.
\end{enumerate}
\end{proposition}

Each condition can be checked easily.
If $\mathbf{x^*_i}$ is an all-or-nothing assignment to shortest paths, then $\mathbf{t(\mathbf{x_i})}^T (\mathbf{x^*_i} - \mathbf{x_i}) \leq 0$, with equality only if $\mathbf{x_i}$ is an equilibrium.
Since $\nabla F(\mathbf{x}) = \mathbf{t}$, and since equilibria correspond to minima of $F$, the first condition is satisfied.
Since $X$ is compact, continuity of the link performance functions ensures that the elements of $HF$ and that the magnitudes of $\mathbf{x^*_i - \mathbf{x}}$ are bounded, and the second condition is satisfied as well.

As a convex program, monotone S-TAP can also be solved by the Frank-Wolfe algorithm, which selects each $\lambda_i \in [0, 1]$ to minimize the value of the objective.
A proof of convergence can be found in Section 2.2.2 of~\cite{bertsekas_nlp}.

\subsection{Algorithms equilibrating paths}

In many path- and bush-based algorithms, the fundamental operation involves equilibrating two paths: given a lower-cost path $\pi_L$ and a higher-cost path $\pi_U$ connecting the same origin and destination, shifting flow from $\pi_U$ to $\pi_L$ to either (approximately) equalize their costs, or to shift all flow onto $\pi_L$.
Examples of such algorithms are the gradient projection (GP) method of~\cite{jaya}, Algorithm B~\citep{Dial}, and TAPAS~\citep{bar2010}.\footnote{TAPAS actually performs such shifts for multiple origin-destination pairs simultaneously, with alternative paths differing on the same segments of links; this point does not affect the discussion here.}

Given these paths, the question is how much flow $\Delta x$ to shift from $\pi_U$ to $\pi_L$.
For TAP, Newton's method is commonly used to estimate the value of $\Delta x$ that equalizes the path costs; this is also the value of $\Delta x$ minimizing the Beckmann function.
The same applies for S-TAP, although the scaling factor in the denominator must change to reflect link interactions, as we now show.

Let $\mathbf{x}(\Delta x)$ denote the link flows after $\Delta x$ has been shifted away from $\pi_U$ onto $\pi_L$.
The only links whose flows will change are those in $\pi_L$ or $\pi_U$, but not both; let $A_L$ and $A_U$ respectively denote the links only in $\pi_L$ and $\pi_U$.
Then
\begin{equation}
\label{deltax}
x_a(\Delta x) = 
\begin{cases}
    x_a + \Delta x & \mbox{if } a \in A_L \\
    x_a - \Delta x & \mbox{if } a \in A_U \\     
    x_a & \mbox{otherwise}
\end{cases}
\,.    
\end{equation}
To find the value of $\Delta x$ minimizing the S-TAP objective $F$, we find where the derivative
\begin{equation}
\label{eq:firstder}
\frac{dF}{d \Delta x} = \sum_a \frac{\partial F}{\partial x_a} \frac{d x_a}{d \Delta x} = \sum_{a \in A_L} t_a - \sum_{a \in A_U} t_a
\end{equation}
vanishes.
As in the separable case, this is exactly when $\pi_L$ and $\pi_U$ have equal cost.

To apply Newton's method, we also need the second derivative of $F$ with respect to $\Delta x$, in order to scale the step size properly. 
We calculate 
\begin{multline}
\frac{d^2F}{d (\Delta x)^2} = \sum_a \sum_{a'} \frac{\partial^2 F}{\partial x_a \partial x_{a'}} \frac{d x_a}{d \Delta x} \frac{d x_{a'}}{d \Delta x} 
\\ = \sum_{a \in A_L} \sum_{a' \in A_L} \frac{\partial t_a}{\partial x_{a'}} + \sum_{a \in A_U} \sum_{a' \in A_U} \frac{\partial t_a}{\partial x_{a'}} - \sum_{a \in A_L} \sum_{a' \in A_U} \frac{\partial t_a}{\partial x_{a'}} - \sum_{a \in A_U} \sum_{a' \in A_L} \frac{\partial t_a}{\partial x_{a'}}
\,.
\end{multline}
Using the symmetry condition, this simplifies to
\begin{multline}
\label{eq:secondder}
\frac{d^2F}{d (\Delta x)^2} = \sum_a \sum_{a'} \frac{\partial^2 F}{\partial x_a \partial x_{a'}} \frac{d x_a}{d \Delta x} \frac{d x_{a'}}{d \Delta x}
\\ = \sum_{a \in A_L} \sum_{a' \in A_L} \frac{\partial t_a}{\partial x_{a'}} + \sum_{a \in A_U} \sum_{a' \in A_U} \frac{\partial t_a}{\partial x_{a'}} - 2 \sum_{a \in A_L} \sum_{a' \in A_U} \frac{\partial t_a}{\partial x_{a'}} 
\,.
\end{multline}
If there are no interactions at all, this formula reduces to $\sum_{a \in A_L \cup A_U} t'_a$, as it must.

So, the Newton estimate for the flow shift equalizing path costs is the quotient of equations~\eqref{eq:firstder} and~\eqref{eq:secondder}.  
The denominator is strictly positive, because $HF$ is positive definite by monotonicity, and this ratio is well-defined.
To preserve feasibility, the flow shift is truncated if any flow would become negative.
This corresponds to the case when the longer path becomes unused after the shift.

This operation can be substituted for the flow shift step in gradient projection, Algorithm B, or TAPAS.

\section{Numerical results}
\label{sec:numerical}

This section describes tests of the algorithms described in the previous section.
The key questions concern the computation time needed when link interactions are considered, compared to the separable case.
We focus on this issue since computational efficiency is one of the advantages static models have over dynamic ones, and if S-TAP is to be useful in practice it should maintain this advantage.
This section also considers asymmetric instances as well, using these algorithms as a heuristic.

We measure convergence using the relative gap, defined as the ratio between the total system travel time, and the total travel time of an all-or-nothing assignment to shortest paths (keeping the current link costs).
Using $\kappa_{rs}$ to reflect the shortest path cost between nodes $r$ and $s$, the relative gap is calculated as
\begin{equation}
RG = \frac{\sum_a t_a x_a - \sum_{(r,s) \in Z^2} d_{rs} \kappa_{rs} }{\sum_{(r,s) \in Z^2} d_{rs} \kappa_{rs}}
\,.
\end{equation}
Additional details on the relative gap, its relationship with other gap measures, and with measures of effectiveness such as link flows and aggregate travel times, are discussed in \cite{Patil2020}.

\subsection{Motivation - Toy example}
Consider the example network shown in Figure~\ref{fig:toy_example}.
We will use this network with different cost functions to show the effects of symmetric vs.\ asymmetric interactions, and the impact of how many other links affect the travel time of a given link.
Consider the five cases (and the corresponding link costs) shown in Table~\ref{table:toy_example_cases}.

The convergence results for the first few iterations of gradient projection are given in Table \ref{table:toy_numbers}.
More iterations for this toy example are omitted for brevity.
The symmetric-full interaction scenario achieves a lower relative gap for the same number of iterations, followed by asymmetric-full interaction, symmetric-partial interaction, asymmetric-partial interaction, and lastly, separable scenario.

This is due to two effects.
First, a higher level of interaction leads to faster convergence.
We speculate that this occurs because each path flow equilibration affects more links in the network than just those in $A_L$ and $A_U$, and therefore moves more of the network towards equilibrium at each step.
Compare the symmetric-full and asymmetric-full scenarios to the symmetric-partial and asymmetric-partial scenarios in Table \ref{table:toy_numbers}.
The full-interaction cases perform better than partial interaction scenarios, which still perform better than the separable TAP scenario.

Second, symmetric scenarios tend to achieve (somewhat) faster convergence compared to asymmetric scenarios.
This is attributed to the path equilibration step having accurate information about the rest of the network (by symmetric Jacobian effects) as opposed to approximate information in the asymmetric case.
Consider the symmetric-partial and asymmetric-partial cases, with very minor link cost differences on links 2 and 3.
The symmetric-partial relative gap is almost half of the asymmetric-partial gap by iteration 5, a trend that continues in additional iterations not shown in the detail.

These two effects can be quantified using the eigenvalues of the Hessian of the cost matrix. For our problem, this is represented by the weight matrix for linear cost functions, and approximated by it otherwise. For steepest descent methods, the condition number (ratio of the largest eigenvalue to the smallest eigenvalue) correlates to the rate of convergence \citep{bertsekas_nlp}. A large condition number means the problem is ill-conditioned, i.e., the optimization variables are not relatively scaled well, and convergence will be slow. A small condition number (closer to 1) will have faster convergence behavior. For example, the condition numbers for the symmetric and asymmetric toy examples with full interaction in Table \ref{table:toy_example_cases} are 3 and 3.154, respectively.

In our investigations, we observe that problem instances with differing condition numbers follow this behavior, and the cost matrix can be useful for predicting convergence behavior w.r.t. similar problem instances. For instance, see Figure \ref{fig:condition_num}. We generated problem instances with slightly different weight matrices (and therefore, condition numbers), and then allowed GP to solve the instances for 20 iterations. We can see that the instances with lower condition numbers generally show better convergence than instances with higher condition numbers.

The next subsection details our experiment design on significantly larger real-world networks to corroborate these observations, and our intuition about convergence for the symmetric/asymmetric and no interaction/full interaction cases.

\subsection{Data and Experiment design}
\label{sec:exp_design}

We test the method of successive averages (MSA), Frank-Wolfe (FW), and gradient projection (GP) on five standard networks, chosen for their varying size and congestion levels.
The networks are obtained from the transportation problems test repository~\citep{bstabler}. Our implementations of MSA, FW, and GP for S-TAP and A-TAP, as well as the testing framework, can be found on the first author's Github repository~\citep{github_patil}.
Table \ref{table:network_data} contains the network size details and average link volume over capacity as a proxy for congestion.
These experiments used a relative gap of $10^{-6}$ as a convergence criterion.

We used the following procedure to generate cost functions with interactions that attempt to preserve the level of congestion in the original networks.
The link performance functions in the original network are separable.
In our experiments, we replace each $x_a$ with a linear combination of the form $\sum_{a' \in A} w_{aa'} x_{a'}$, with $\sum_{a'} w_{aa'} = 1$ and each $w_{aa'} \geq 0$.  
If the weight matrix $\mathbf{W}$ is symmetric, then the interactions are approximately symmetric (but not entirely so, since the link performance functions are nonlinear).
The separable case is represented with $\mathbf{W} = \mathbf{I}$.
We generate the weight matrices so that each link depends on a given number of other links (the number of ``degrees of dependency,'' denoted $N$), and with $\mathbf{W}$ diagonally dominant to avoid cases with multiple equilibria.

The first set of experiments tested the convergence behavior of the three algorithms for TAP, S-TAP, and A-TAP on all networks.
The best performing algorithm was then chosen for further convergence testing of S-TAP and A-TAP.
The second set of experiments studied the effect of topographical link interactions.
These experiments aim to understand impact the degree of dependency has on the convergence rate.
The third set of experiments analyzed the effects of symmetry.  
We consider a smooth transition from asymmetry to symmetry to starting with an  asymmetric matrix $\textbf{W}$, and taking weighted averages with the associated symmetric matrix $\frac{1}{2}(\mathbf{W} + \mathbf{W}^T)$.
Specifically, these experiments consider the Jacobian matrices $\lambda \textbf{W} + (1-\lambda)\frac{1}{2}(\mathbf{W} + \mathbf{W}^T)$ for $\lambda \in \{0, 0.2, 0.4, 0.6, 0.8, 1 \}$.

\subsection{Results}

Figures \ref{fig:SF_compiled_results} and \ref{fig:SF_asym_to_sym} present the compiled results for Sioux Falls.
The first observation from Figure \ref{fig:SF_a} is the behavior of MSA and FW algorithms is extremely similar for TAP, S-TAP, and A-TAP, while GP outperforms them.
The relative gap and number of iterations are linearly related on a logarithmic axis, consistent with prior literature~\citep{xie2018}.
Based on these initial results, we use GP as the testing algorithm for further experiments.

An important observation here is the independence from implementation details and absolute computation time.
Our experiments have been conducted on a basic GP implementation.
Therefore, any performance gains achieved by parallelization or other implementation techniques are applicable to these results, helping speed up convergence.
For instance, \cite{chen2020parallel} implement a parallel block coordinate descent algorithm based on GP algorithm, and the absolute computation time gains would be applicable to our observations for S-TAP and A-TAP.
Therefore, we do not emphasize absolute computation times, but focus on comparative behavior. 

Figures \ref{fig:SF_b} and \ref{fig:SF_c} show the behavior of link cost dependency on link flows within $N$ degrees of any given link for S-TAP and A-TAP, respectively.
A higher degree of dependence leads to earlier convergence, as was observed in the motivating toy network.
This is attributed to each flow equilibration step having more implicit information about the network state.
Also, S-TAP is seen to converge marginally faster than A-TAP, as expected.
The only exception is the case $N = 2$ case, though it still does not outperform $N=5$ or $N=6$ link interaction instances.
Along similar lines, Figure \ref{fig:SF_asym_to_sym} shows convergence behavior when the weight matrix proceeds from an asymmetric instance ($\lambda = 1$) to the corresponding symmetric instance ($\lambda = 0$).
The results indicate as the matrix tends toward the symmetric version, more accurate information (about the remaining network) is available during each flow equilibration, leading to marginally faster convergence.

Figures \ref{fig:EMA_compiled_results}--\ref{fig:CR_compiled_results} show the results for Eastern Massachusetts, Chicago Sketch, Barcelona, and Chicago Regional networks, respectively.
The main observations from the Sioux Falls networks are seen to hold for these larger networks.
The maximum $N$ for these four networks is 8, 31, 30, and 112, respectively.
These networks show no exception to these trends, unlike the Sioux Falls network.
Also, the test instances with lower $N$ values achieved better relative gap levels for larger networks; the highest $N$ instances all reached a gap below $10^{-12}$.

\section{Convergence experiments}
\label{sec:algconv_conv_exps}

This subsection focuses on the convergence behavior for S-TAP. This is particularly important for applications where TAP is solved multiple times as a sub-problem. Some applications include network design problems~\citep{yang, patil2022budget, gokalp2021post}, sensitivity analysis~\citep{boyles3}, trip table estimation~\citep{yang1995heuristic, lundgren2008heuristic}, network pricing~\citep{yan1996optimal}, and other bilevel optimization problems \citep{yin2000genetic, astroza, pandey2022computationally, boyles}. \cite{Patil2020} provide additional background on convergence analysis for TAP. 

Here, we primarily consider GP with two-link interactions. This limits the scope of our experiments to the most common use case, i.e., two-way roads. The networks studied in this paper are shown in Table \ref{table:nets_STAP}, all obtained from the Transportation Networks for Research repository \citep{bstabler}. For ease of reference, we categorize the networks roughly by size: Sioux Falls through Anaheim are designated as small, Chicago Sketch through Terrassa are designated as medium, and the remaining networks are designated as large.  The last column in this table shows the average equilibrium flow-to-capacity ratios, excluding centroid connectors.  We consider networks with ratios of less than 0.5 to be uncongested, with ratios between 0.5 and 1.0 to be semi-congested, and networks with ratios greater than 1.0 to be congested.  The Terrassa network is a clear outlier in this regard, assigning over 25 million trips in a region whose current population is around 200,000, resulting in a flow-to-capacity ratio of almost 6.  While such a demand level may not be realistic, we nevertheless include this network as a ``stress test'' to see whether consistent trends can be seen even in extremely congested networks. Lastly, we also conduct these experiments using Algorithm B (AlgB), a bush-based algorithm which is often faster in practice.

Given a feasible solution $(\mathbf{x}, \mathbf{h})$ to TAP, we select three metrics for analysis.  The \emph{total system travel time} (TSTT) expresses the sum of each vehicle's travel time in the network as

\begin{equation}
    TSTT(\mathbf{x}) = \sum_{(i,j) \in A} t_{ij}x_{ij}\,.\\
\end{equation}

\emph{Vehicle-miles traveled} (VMT) expresses the total distance traveled by vehicles in the network as

\begin{equation}
    VMT(\mathbf{x}) = \sum_{(i,j) \in A} l_{ij}x_{ij}\,.\\
\end{equation}

To measure convergence of these metrics, we calculate the relative difference between their values at the current solution $\mathbf{x}$ and the equilibrium solution $\mathbf{x^*}$ as:

\begin{equation}
    \Delta TSTT(\mathbf{x}) = \frac{TSTT(\mathbf{x}) - TSTT(\mathbf{x^*})}{TSTT(\mathbf{x^*})}
\end{equation}

\begin{equation}
    \Delta VMT(\mathbf{x}) = \frac{VMT(\mathbf{x}) - VMT(\mathbf{x^*})}{VMT(\mathbf{x^*})}  \,.
\end{equation}

Both TSTT and VMT are aggregate metrics.  To represent convergence of the specific link and path flows themselves, we measure the proportion of links within a given relative threshold $\epsilon$ of their equilibrium values.  Let $A_\epsilon^*(\mathbf{x})$ denote the set of links with flows within this threshold:

\begin{equation}
    A^*_\epsilon(\mathbf{x}) = \left\{ (i,j) \in A : \left| x_{ij} - x^*_{ij} \right| < \epsilon x^*_{ij} \right\} \,.
\end{equation}

Using these sets, we define the \emph{proportion of unconverged links (PUL)} as

\begin{equation}
    PUL(\mathbf{x},\epsilon) = 1 - \frac{|A^*_\epsilon(\mathbf{x})|}{|A|} \,.
\end{equation}

These metrics --- $\Delta TSTT$, $\Delta VMT$, $PUL$ --- are directly related and used in practical applications and planning, and converge to zero at the equilibrium solution. We track these metrics against relative gap (defined in Section \ref{sec:numerical}), a convergence metric.

The full set of results can be found in the Appendix (Tables \ref{tab:STAP-gp}-\ref{tab:algB_GP}). The values of the three main metrics for GP experiments are shown in Figures \ref{fig:TSTT_trend}--\ref{fig:PUL_trend} (presented according to each metric).  The thin lines represent the values of each metric in one of the twelve networks tested, and the thick line represents the average value.  Both sets of figures use logarithmic axes for the relative gap. All metrics converged at roughly similar rates, despite significant differences in the size and congestion level of the networks tested.

In all the networks, the aggregate metrics ($TSTT$ and $VMT$) are already very near stabilization at a relative gap of $10^{-3}$.  For the small and medium networks, these values are within 1\% of the equilibrium values when the relative gap is $10^{-4}$, and for the large networks they are within 2\%.  Both $\Delta TSTT$ and $\Delta VMT$ converge at roughly similar rates, but $\Delta VMT$ is usually slightly lower at a particular gap level. This behavior is in line with the metric behavior for TAP, as noted in \cite{Patil2020}.

The proportion of unconverged links was the metric originally studied by \cite{boyce} for the Philadelphia regional network.  They found that a gap of $10^{-4}$ was required to approach convergence for freeway links, defining convergence as a $PUL$ of 1\% or less.  To achieve this level of convergence for arterial links as well as freeway links, a relative gap of $10^{-5}$ was needed.  Our results show that this latter conclusion generally holds across the other networks tested, and that 99\% of link flows are accurate to within 1\% of equilibrium values at this gap level.

Next, we compare the GP convergence behavior to AlgB convergence behavior. The full data from these results are shown in Tables~\ref{tab:STAP-algb} (raw data for Algorithm B) and~\ref{tab:algB_GP} (for a side-by-side comparison) in the Appendix. The trends are very similar between the two algorithms, and the values of each metric are always of the same order of magnitude, and almost always nearly identical numerically.  This finding is encouraging, suggesting that the conclusions of the GP experiments are applicable to other algorithms, and that the relative gap is a good universal measure of convergence, regardless of the specific assignment algorithm. Again, these results are in line with prior results from \cite{Patil2020}.

\section{Conclusions}
\label{sec:algconv_conclusions}

This study reconsiders the traffic assignment problems with interactions.
We showed that merge models from dynamic traffic assignment can be approximated with symmetric, monotone cost functions; examined alternative algorithms, deriving the Newton shift formula with interactions; and considered the practical convergence rate of gradient projection and Algorithm B with such functions.
We found that instances with interactions actually converge faster than separable instances of traffic assignment.
All of this suggests that problems currently studied with static traffic assignment may benefit from considering interactions --- there seems to be little computational difficulty (in fact, convergence was almost always faster), and the valuable properties of equilibrium existence and uniqueness are retained.
Of course, there remain other problems where dynamic assignment is preferred.

Further research would be valuable along several lines.
Further investigation of appropriate cost functions is needed, to derive them from other node models, and to consider the impacts of these approximations on network-wide flow.
Additional research into algorithms for the non-monotone or non-symmetric cases is also needed; we show that gradient projection still functions as an acceptable heuristic, but comparisons with other algorithms from the variational inequality literature are needed.

\section*{Acknowledgments}

This research was supported by the National Science Foundation under Grant CMMI-1826320.

\section*{Disclosure statement}

 The authors report there are no competing interests to declare.

\section*{Data availability statement}

The authors confirm that the data supporting the findings of this study are available within the article and its supplementary materials.

\clearpage

\bibliographystyle{tfcad}
\bibliography{interactcadsample}
\newpage

\appendix
\section{Raw Data}

\begin{landscape}
\begin{table}[]
\centering
\caption{Metric stabilization behavior data using Gradient Projection for S-TAP}
\label{tab:STAP-gp}
\small
\begin{tabular}{cccccccccccccc}
\cline{1-4} \cline{6-9} \cline{11-14}
\multicolumn{4}{|c|}{\textbf{Sioux Falls}} & \multicolumn{1}{c|}{} & \multicolumn{4}{c|}{\textbf{Berlin-Mitte-Prenzlauerberg-Friedrichshain-Center}} & \multicolumn{1}{c|}{} & \multicolumn{4}{c|}{\textbf{Austin}} \\ \cline{1-4} \cline{6-9} \cline{11-14} 
\multicolumn{1}{|c|}{\textbf{Gap Level}} & \multicolumn{1}{c|}{\textbf{$\Delta$TSTT}} & \multicolumn{1}{c|}{\textbf{$\Delta$VMT}} & \multicolumn{1}{c|}{\textbf{PUL}} & \multicolumn{1}{c|}{} & \multicolumn{1}{c|}{\textbf{Gap Level}} & \multicolumn{1}{c|}{\textbf{$\Delta$TSTT}} & \multicolumn{1}{c|}{\textbf{$\Delta$VMT}} & \multicolumn{1}{c|}{\textbf{PUL}} & \multicolumn{1}{c|}{} & \multicolumn{1}{c|}{\textbf{Gap Level}} & \multicolumn{1}{c|}{\textbf{$\Delta$TSTT}} & \multicolumn{1}{c|}{\textbf{$\Delta$VMT}} & \multicolumn{1}{c|}{\textbf{PUL}} \\ \cline{1-4} \cline{6-9} \cline{11-14} 
\multicolumn{1}{|c|}{1E-03} & \multicolumn{1}{c|}{2.13\%} & \multicolumn{1}{c|}{1.48\%} & \multicolumn{1}{c|}{5.26\%} & \multicolumn{1}{c|}{} & \multicolumn{1}{c|}{1E-03} & \multicolumn{1}{c|}{2.04\%} & \multicolumn{1}{c|}{1.42\%} & \multicolumn{1}{c|}{6.24\%} & \multicolumn{1}{c|}{} & \multicolumn{1}{c|}{1E-03} & \multicolumn{1}{c|}{2.39\%} & \multicolumn{1}{c|}{2.18\%} & \multicolumn{1}{c|}{9.55\%} \\ \cline{1-4} \cline{6-9} \cline{11-14} 
\multicolumn{1}{|c|}{1E-04} & \multicolumn{1}{c|}{0.89\%} & \multicolumn{1}{c|}{0.65\%} & \multicolumn{1}{c|}{2.63\%} & \multicolumn{1}{c|}{} & \multicolumn{1}{c|}{1E-04} & \multicolumn{1}{c|}{0.98\%} & \multicolumn{1}{c|}{0.58\%} & \multicolumn{1}{c|}{1.95\%} & \multicolumn{1}{c|}{} & \multicolumn{1}{c|}{1E-04} & \multicolumn{1}{c|}{1.59\%} & \multicolumn{1}{c|}{1.09\%} & \multicolumn{1}{c|}{5.38\%} \\ \cline{1-4} \cline{6-9} \cline{11-14} 
\multicolumn{1}{|c|}{1E-05} & \multicolumn{1}{c|}{0.35\%} & \multicolumn{1}{c|}{0.18\%} & \multicolumn{1}{c|}{1.31\%} & \multicolumn{1}{c|}{} & \multicolumn{1}{c|}{1E-05} & \multicolumn{1}{c|}{0.36\%} & \multicolumn{1}{c|}{0.23\%} & \multicolumn{1}{c|}{0.78\%} & \multicolumn{1}{c|}{} & \multicolumn{1}{c|}{1E-05} & \multicolumn{1}{c|}{0.85\%} & \multicolumn{1}{c|}{0.59\%} & \multicolumn{1}{c|}{0.90\%} \\ \cline{1-4} \cline{6-9} \cline{11-14} 
\multicolumn{1}{|c|}{1E-06} & \multicolumn{1}{c|}{0.07\%} & \multicolumn{1}{c|}{0.03\%} & \multicolumn{1}{c|}{0.00\%} & \multicolumn{1}{c|}{} & \multicolumn{1}{c|}{1E-06} & \multicolumn{1}{c|}{0.23\%} & \multicolumn{1}{c|}{0.17\%} & \multicolumn{1}{c|}{0.44\%} & \multicolumn{1}{c|}{} & \multicolumn{1}{c|}{1E-06} & \multicolumn{1}{c|}{0.58\%} & \multicolumn{1}{c|}{0.39\%} & \multicolumn{1}{c|}{0.72\%} \\ \cline{1-4} \cline{6-9} \cline{11-14} 
\multicolumn{1}{|c|}{1E-07} & \multicolumn{1}{c|}{0.02\%} & \multicolumn{1}{c|}{0.01\%} & \multicolumn{1}{c|}{0.00\%} & \multicolumn{1}{c|}{} & \multicolumn{1}{c|}{1E-07} & \multicolumn{1}{c|}{0.09\%} & \multicolumn{1}{c|}{0.08\%} & \multicolumn{1}{c|}{0.10\%} & \multicolumn{1}{c|}{} & \multicolumn{1}{c|}{1E-07} & \multicolumn{1}{c|}{0.19\%} & \multicolumn{1}{c|}{0.19\%} & \multicolumn{1}{c|}{0.19\%} \\ \cline{1-4} \cline{6-9} \cline{11-14} 
\multicolumn{1}{|c|}{1E-08} & \multicolumn{1}{c|}{0.00\%} & \multicolumn{1}{c|}{0.00\%} & \multicolumn{1}{c|}{0.00\%} & \multicolumn{1}{c|}{} & \multicolumn{1}{c|}{1E-08} & \multicolumn{1}{c|}{0.05\%} & \multicolumn{1}{c|}{0.02\%} & \multicolumn{1}{c|}{0.05\%} & \multicolumn{1}{c|}{} & \multicolumn{1}{c|}{1E-08} & \multicolumn{1}{c|}{0.05\%} & \multicolumn{1}{c|}{0.03\%} & \multicolumn{1}{c|}{0.08\%} \\ \cline{1-4} \cline{6-9} \cline{11-14} 
 &  &  &  &  &  &  &  &  &  &  &  &  &  \\ \cline{1-4} \cline{6-9} \cline{11-14} 
\multicolumn{4}{|c|}{\textbf{Eastern Massachusetts}} & \multicolumn{1}{c|}{} & \multicolumn{4}{c|}{\textbf{Barcelona}} & \multicolumn{1}{c|}{} & \multicolumn{4}{c|}{\textbf{Berlin Center}} \\ \cline{1-4} \cline{6-9} \cline{11-14} 
\multicolumn{1}{|c|}{\textbf{Gap Level}} & \multicolumn{1}{c|}{\textbf{$\Delta$TSTT}} & \multicolumn{1}{c|}{\textbf{$\Delta$VMT}} & \multicolumn{1}{c|}{\textbf{PUL}} & \multicolumn{1}{c|}{} & \multicolumn{1}{c|}{\textbf{Gap Level}} & \multicolumn{1}{c|}{\textbf{$\Delta$TSTT}} & \multicolumn{1}{c|}{\textbf{$\Delta$VMT}} & \multicolumn{1}{c|}{\textbf{PUL}} & \multicolumn{1}{c|}{} & \multicolumn{1}{c|}{\textbf{Gap Level}} & \multicolumn{1}{c|}{\textbf{$\Delta$TSTT}} & \multicolumn{1}{c|}{\textbf{$\Delta$VMT}} & \multicolumn{1}{c|}{\textbf{PUL}} \\ \cline{1-4} \cline{6-9} \cline{11-14} 
\multicolumn{1}{|c|}{1E-03} & \multicolumn{1}{c|}{1.80\%} & \multicolumn{1}{c|}{1.25\%} & \multicolumn{1}{c|}{4.50\%} & \multicolumn{1}{c|}{} & \multicolumn{1}{c|}{1E-03} & \multicolumn{1}{c|}{1.78\%} & \multicolumn{1}{c|}{1.65\%} & \multicolumn{1}{c|}{5.52\%} & \multicolumn{1}{c|}{} & \multicolumn{1}{c|}{1E-03} & \multicolumn{1}{c|}{2.58\%} & \multicolumn{1}{c|}{2.29\%} & \multicolumn{1}{c|}{10.54\%} \\ \cline{1-4} \cline{6-9} \cline{11-14} 
\multicolumn{1}{|c|}{1E-04} & \multicolumn{1}{c|}{0.75\%} & \multicolumn{1}{c|}{0.44\%} & \multicolumn{1}{c|}{2.66\%} & \multicolumn{1}{c|}{} & \multicolumn{1}{c|}{1E-04} & \multicolumn{1}{c|}{0.94\%} & \multicolumn{1}{c|}{0.93\%} & \multicolumn{1}{c|}{2.54\%} & \multicolumn{1}{c|}{} & \multicolumn{1}{c|}{1E-04} & \multicolumn{1}{c|}{1.87\%} & \multicolumn{1}{c|}{1.13\%} & \multicolumn{1}{c|}{4.61\%} \\ \cline{1-4} \cline{6-9} \cline{11-14} 
\multicolumn{1}{|c|}{1E-05} & \multicolumn{1}{c|}{0.29\%} & \multicolumn{1}{c|}{0.10\%} & \multicolumn{1}{c|}{0.87\%} & \multicolumn{1}{c|}{} & \multicolumn{1}{c|}{1E-05} & \multicolumn{1}{c|}{0.45\%} & \multicolumn{1}{c|}{0.43\%} & \multicolumn{1}{c|}{0.91\%} & \multicolumn{1}{c|}{} & \multicolumn{1}{c|}{1E-05} & \multicolumn{1}{c|}{0.99\%} & \multicolumn{1}{c|}{0.65\%} & \multicolumn{1}{c|}{1.20\%} \\ \cline{1-4} \cline{6-9} \cline{11-14} 
\multicolumn{1}{|c|}{1E-06} & \multicolumn{1}{c|}{0.20\%} & \multicolumn{1}{c|}{0.06\%} & \multicolumn{1}{c|}{0.72\%} & \multicolumn{1}{c|}{} & \multicolumn{1}{c|}{1E-06} & \multicolumn{1}{c|}{0.23\%} & \multicolumn{1}{c|}{0.31\%} & \multicolumn{1}{c|}{0.67\%} & \multicolumn{1}{c|}{} & \multicolumn{1}{c|}{1E-06} & \multicolumn{1}{c|}{0.64\%} & \multicolumn{1}{c|}{0.44\%} & \multicolumn{1}{c|}{0.79\%} \\ \cline{1-4} \cline{6-9} \cline{11-14} 
\multicolumn{1}{|c|}{1E-07} & \multicolumn{1}{c|}{0.05\%} & \multicolumn{1}{c|}{0.03\%} & \multicolumn{1}{c|}{0.41\%} & \multicolumn{1}{c|}{} & \multicolumn{1}{c|}{1E-07} & \multicolumn{1}{c|}{0.12\%} & \multicolumn{1}{c|}{0.11\%} & \multicolumn{1}{c|}{0.29\%} & \multicolumn{1}{c|}{} & \multicolumn{1}{c|}{1E-07} & \multicolumn{1}{c|}{0.19\%} & \multicolumn{1}{c|}{0.16\%} & \multicolumn{1}{c|}{0.33\%} \\ \cline{1-4} \cline{6-9} \cline{11-14} 
\multicolumn{1}{|c|}{1E-08} & \multicolumn{1}{c|}{0.03\%} & \multicolumn{1}{c|}{0.01\%} & \multicolumn{1}{c|}{0.04\%} & \multicolumn{1}{c|}{} & \multicolumn{1}{c|}{1E-08} & \multicolumn{1}{c|}{0.02\%} & \multicolumn{1}{c|}{0.03\%} & \multicolumn{1}{c|}{0.08\%} & \multicolumn{1}{c|}{} & \multicolumn{1}{c|}{1E-08} & \multicolumn{1}{c|}{0.08\%} & \multicolumn{1}{c|}{0.03\%} & \multicolumn{1}{c|}{0.08\%} \\ \cline{1-4} \cline{6-9} \cline{11-14} 
 &  &  &  &  &  &  &  &  &  &  &  &  &  \\ \cline{1-4} \cline{6-9} \cline{11-14} 
\multicolumn{4}{|c|}{\textbf{Anaheim}} & \multicolumn{1}{c|}{} & \multicolumn{4}{c|}{\textbf{Winnipeg}} & \multicolumn{1}{c|}{} & \multicolumn{4}{c|}{\textbf{Chicago-Regional}} \\ \cline{1-4} \cline{6-9} \cline{11-14} 
\multicolumn{1}{|c|}{\textbf{Gap Level}} & \multicolumn{1}{c|}{\textbf{$\Delta$TSTT}} & \multicolumn{1}{c|}{\textbf{$\Delta$VMT}} & \multicolumn{1}{c|}{\textbf{PUL}} & \multicolumn{1}{c|}{} & \multicolumn{1}{c|}{\textbf{Gap Level}} & \multicolumn{1}{c|}{\textbf{$\Delta$TSTT}} & \multicolumn{1}{c|}{\textbf{$\Delta$VMT}} & \multicolumn{1}{c|}{\textbf{PUL}} & \multicolumn{1}{c|}{} & \multicolumn{1}{c|}{\textbf{Gap Level}} & \multicolumn{1}{c|}{\textbf{$\Delta$TSTT}} & \multicolumn{1}{c|}{\textbf{$\Delta$VMT}} & \multicolumn{1}{c|}{\textbf{PUL}} \\ \cline{1-4} \cline{6-9} \cline{11-14} 
\multicolumn{1}{|c|}{1E-03} & \multicolumn{1}{c|}{2.47\%} & \multicolumn{1}{c|}{1.41\%} & \multicolumn{1}{c|}{4.77\%} & \multicolumn{1}{c|}{} & \multicolumn{1}{c|}{1E-03} & \multicolumn{1}{c|}{2.84\%} & \multicolumn{1}{c|}{1.81\%} & \multicolumn{1}{c|}{5.49\%} & \multicolumn{1}{c|}{} & \multicolumn{1}{c|}{1E-03} & \multicolumn{1}{c|}{2.24\%} & \multicolumn{1}{c|}{1.51\%} & \multicolumn{1}{c|}{11.58\%} \\ \cline{1-4} \cline{6-9} \cline{11-14} 
\multicolumn{1}{|c|}{1E-04} & \multicolumn{1}{c|}{0.77\%} & \multicolumn{1}{c|}{0.57\%} & \multicolumn{1}{c|}{2.58\%} & \multicolumn{1}{c|}{} & \multicolumn{1}{c|}{1E-04} & \multicolumn{1}{c|}{0.94\%} & \multicolumn{1}{c|}{1.01\%} & \multicolumn{1}{c|}{2.20\%} & \multicolumn{1}{c|}{} & \multicolumn{1}{c|}{1E-04} & \multicolumn{1}{c|}{1.24\%} & \multicolumn{1}{c|}{0.84\%} & \multicolumn{1}{c|}{5.07\%} \\ \cline{1-4} \cline{6-9} \cline{11-14} 
\multicolumn{1}{|c|}{1E-05} & \multicolumn{1}{c|}{0.41\%} & \multicolumn{1}{c|}{0.31\%} & \multicolumn{1}{c|}{0.72\%} & \multicolumn{1}{c|}{} & \multicolumn{1}{c|}{1E-05} & \multicolumn{1}{c|}{0.47\%} & \multicolumn{1}{c|}{0.30\%} & \multicolumn{1}{c|}{0.91\%} & \multicolumn{1}{c|}{} & \multicolumn{1}{c|}{1E-05} & \multicolumn{1}{c|}{0.61\%} & \multicolumn{1}{c|}{0.43\%} & \multicolumn{1}{c|}{0.94\%} \\ \cline{1-4} \cline{6-9} \cline{11-14} 
\multicolumn{1}{|c|}{1E-06} & \multicolumn{1}{c|}{0.25\%} & \multicolumn{1}{c|}{0.18\%} & \multicolumn{1}{c|}{0.24\%} & \multicolumn{1}{c|}{} & \multicolumn{1}{c|}{1E-06} & \multicolumn{1}{c|}{0.26\%} & \multicolumn{1}{c|}{0.30\%} & \multicolumn{1}{c|}{0.52\%} & \multicolumn{1}{c|}{} & \multicolumn{1}{c|}{1E-06} & \multicolumn{1}{c|}{0.46\%} & \multicolumn{1}{c|}{0.34\%} & \multicolumn{1}{c|}{0.77\%} \\ \cline{1-4} \cline{6-9} \cline{11-14} 
\multicolumn{1}{|c|}{1E-07} & \multicolumn{1}{c|}{0.10\%} & \multicolumn{1}{c|}{0.08\%} & \multicolumn{1}{c|}{0.04\%} & \multicolumn{1}{c|}{} & \multicolumn{1}{c|}{1E-07} & \multicolumn{1}{c|}{0.11\%} & \multicolumn{1}{c|}{0.10\%} & \multicolumn{1}{c|}{0.16\%} & \multicolumn{1}{c|}{} & \multicolumn{1}{c|}{1E-07} & \multicolumn{1}{c|}{0.18\%} & \multicolumn{1}{c|}{0.08\%} & \multicolumn{1}{c|}{0.27\%} \\ \cline{1-4} \cline{6-9} \cline{11-14} 
\multicolumn{1}{|c|}{1E-08} & \multicolumn{1}{c|}{0.05\%} & \multicolumn{1}{c|}{0.01\%} & \multicolumn{1}{c|}{0.02\%} & \multicolumn{1}{c|}{} & \multicolumn{1}{c|}{1E-08} & \multicolumn{1}{c|}{0.02\%} & \multicolumn{1}{c|}{0.03\%} & \multicolumn{1}{c|}{0.05\%} & \multicolumn{1}{c|}{} & \multicolumn{1}{c|}{1E-08} & \multicolumn{1}{c|}{0.04\%} & \multicolumn{1}{c|}{0.02\%} & \multicolumn{1}{c|}{0.04\%} \\ \cline{1-4} \cline{6-9} \cline{11-14} 
 &  &  &  &  &  &  &  &  &  &  &  &  &  \\ \cline{1-4} \cline{6-9} \cline{11-14} 
\multicolumn{4}{|c|}{\textbf{Chicago Sketch}} & \multicolumn{1}{c|}{} & \multicolumn{4}{c|}{\textbf{Terrassa}} & \multicolumn{1}{c|}{} & \multicolumn{4}{c|}{\textbf{Philadelphia}} \\ \cline{1-4} \cline{6-9} \cline{11-14} 
\multicolumn{1}{|c|}{\textbf{Gap Level}} & \multicolumn{1}{c|}{\textbf{$\Delta$TSTT}} & \multicolumn{1}{c|}{\textbf{$\Delta$VMT}} & \multicolumn{1}{c|}{\textbf{PUL}} & \multicolumn{1}{c|}{} & \multicolumn{1}{c|}{\textbf{Gap Level}} & \multicolumn{1}{c|}{\textbf{$\Delta$TSTT}} & \multicolumn{1}{c|}{\textbf{$\Delta$VMT}} & \multicolumn{1}{c|}{\textbf{PUL}} & \multicolumn{1}{c|}{} & \multicolumn{1}{c|}{\textbf{Gap Level}} & \multicolumn{1}{c|}{\textbf{$\Delta$TSTT}} & \multicolumn{1}{c|}{\textbf{$\Delta$VMT}} & \multicolumn{1}{c|}{\textbf{PUL}} \\ \cline{1-4} \cline{6-9} \cline{11-14} 
\multicolumn{1}{|c|}{1E-03} & \multicolumn{1}{c|}{1.88\%} & \multicolumn{1}{c|}{0.94\%} & \multicolumn{1}{c|}{4.87\%} & \multicolumn{1}{c|}{} & \multicolumn{1}{c|}{1E-03} & \multicolumn{1}{c|}{2.55\%} & \multicolumn{1}{c|}{2.16\%} & \multicolumn{1}{c|}{6.51\%} & \multicolumn{1}{c|}{} & \multicolumn{1}{c|}{1E-03} & \multicolumn{1}{c|}{3.83\%} & \multicolumn{1}{c|}{2.43\%} & \multicolumn{1}{c|}{15.82\%} \\ \cline{1-4} \cline{6-9} \cline{11-14} 
\multicolumn{1}{|c|}{1E-04} & \multicolumn{1}{c|}{0.89\%} & \multicolumn{1}{c|}{0.45\%} & \multicolumn{1}{c|}{2.49\%} & \multicolumn{1}{c|}{} & \multicolumn{1}{c|}{1E-04} & \multicolumn{1}{c|}{0.98\%} & \multicolumn{1}{c|}{0.96\%} & \multicolumn{1}{c|}{2.77\%} & \multicolumn{1}{c|}{} & \multicolumn{1}{c|}{1E-04} & \multicolumn{1}{c|}{1.62\%} & \multicolumn{1}{c|}{1.24\%} & \multicolumn{1}{c|}{7.53\%} \\ \cline{1-4} \cline{6-9} \cline{11-14} 
\multicolumn{1}{|c|}{1E-05} & \multicolumn{1}{c|}{0.38\%} & \multicolumn{1}{c|}{0.25\%} & \multicolumn{1}{c|}{0.96\%} & \multicolumn{1}{c|}{} & \multicolumn{1}{c|}{1E-05} & \multicolumn{1}{c|}{0.52\%} & \multicolumn{1}{c|}{0.36\%} & \multicolumn{1}{c|}{0.98\%} & \multicolumn{1}{c|}{} & \multicolumn{1}{c|}{1E-05} & \multicolumn{1}{c|}{1.05\%} & \multicolumn{1}{c|}{0.79\%} & \multicolumn{1}{c|}{1.26\%} \\ \cline{1-4} \cline{6-9} \cline{11-14} 
\multicolumn{1}{|c|}{1E-06} & \multicolumn{1}{c|}{0.29\%} & \multicolumn{1}{c|}{0.08\%} & \multicolumn{1}{c|}{0.32\%} & \multicolumn{1}{c|}{} & \multicolumn{1}{c|}{1E-06} & \multicolumn{1}{c|}{0.29\%} & \multicolumn{1}{c|}{0.37\%} & \multicolumn{1}{c|}{0.48\%} & \multicolumn{1}{c|}{} & \multicolumn{1}{c|}{1E-06} & \multicolumn{1}{c|}{0.63\%} & \multicolumn{1}{c|}{0.35\%} & \multicolumn{1}{c|}{0.82\%} \\ \cline{1-4} \cline{6-9} \cline{11-14} 
\multicolumn{1}{|c|}{1E-07} & \multicolumn{1}{c|}{0.08\%} & \multicolumn{1}{c|}{0.04\%} & \multicolumn{1}{c|}{0.18\%} & \multicolumn{1}{c|}{} & \multicolumn{1}{c|}{1E-07} & \multicolumn{1}{c|}{0.13\%} & \multicolumn{1}{c|}{0.08\%} & \multicolumn{1}{c|}{0.20\%} & \multicolumn{1}{c|}{} & \multicolumn{1}{c|}{1E-07} & \multicolumn{1}{c|}{0.36\%} & \multicolumn{1}{c|}{0.14\%} & \multicolumn{1}{c|}{0.37\%} \\ \cline{1-4} \cline{6-9} \cline{11-14} 
\multicolumn{1}{|c|}{1E-08} & \multicolumn{1}{c|}{0.04\%} & \multicolumn{1}{c|}{0.02\%} & \multicolumn{1}{c|}{0.07\%} & \multicolumn{1}{c|}{} & \multicolumn{1}{c|}{1E-08} & \multicolumn{1}{c|}{0.02\%} & \multicolumn{1}{c|}{0.03\%} & \multicolumn{1}{c|}{0.08\%} & \multicolumn{1}{c|}{} & \multicolumn{1}{c|}{1E-08} & \multicolumn{1}{c|}{0.10\%} & \multicolumn{1}{c|}{0.03\%} & \multicolumn{1}{c|}{0.18\%} \\ \cline{1-4} \cline{6-9} \cline{11-14} 
\end{tabular}

\end{table}

\pagebreak

\begin{table}[]
\centering
\caption{Metric stabilization behavior data using Algorithm B for S-TAP}
\label{tab:STAP-algb}
\small
\begin{tabular}{cccccccccccccc}
\cline{1-4} \cline{6-9} \cline{11-14}
\multicolumn{4}{|c|}{\textbf{Sioux Falls}} & \multicolumn{1}{c|}{} & \multicolumn{4}{c|}{\textbf{Berlin-Mitte-Prenzlauerberg-Friedrichshain-Center}} & \multicolumn{1}{c|}{} & \multicolumn{4}{c|}{\textbf{Austin}} \\ \cline{1-4} \cline{6-9} \cline{11-14} 
\multicolumn{1}{|c|}{\textbf{Gap Level}} & \multicolumn{1}{c|}{\textbf{$\Delta$TSTT}} & \multicolumn{1}{c|}{\textbf{$\Delta$VMT}} & \multicolumn{1}{c|}{\textbf{PUL}} & \multicolumn{1}{c|}{} & \multicolumn{1}{c|}{\textbf{Gap Level}} & \multicolumn{1}{c|}{\textbf{$\Delta$TSTT}} & \multicolumn{1}{c|}{\textbf{$\Delta$VMT}} & \multicolumn{1}{c|}{\textbf{PUL}} & \multicolumn{1}{c|}{} & \multicolumn{1}{c|}{\textbf{Gap Level}} & \multicolumn{1}{c|}{\textbf{$\Delta$TSTT}} & \multicolumn{1}{c|}{\textbf{$\Delta$VMT}} & \multicolumn{1}{c|}{\textbf{PUL}} \\ \cline{1-4} \cline{6-9} \cline{11-14} 
\multicolumn{1}{|c|}{1E-03} & \multicolumn{1}{c|}{2.00\%} & \multicolumn{1}{c|}{1.44\%} & \multicolumn{1}{c|}{5.26\%} & \multicolumn{1}{c|}{} & \multicolumn{1}{c|}{1E-03} & \multicolumn{1}{c|}{2.11\%} & \multicolumn{1}{c|}{1.56\%} & \multicolumn{1}{c|}{6.90\%} & \multicolumn{1}{c|}{} & \multicolumn{1}{c|}{1E-03} & \multicolumn{1}{c|}{2.87\%} & \multicolumn{1}{c|}{2.47\%} & \multicolumn{1}{c|}{11.24\%} \\ \cline{1-4} \cline{6-9} \cline{11-14} 
\multicolumn{1}{|c|}{1E-04} & \multicolumn{1}{c|}{0.80\%} & \multicolumn{1}{c|}{0.59\%} & \multicolumn{1}{c|}{2.63\%} & \multicolumn{1}{c|}{} & \multicolumn{1}{c|}{1E-04} & \multicolumn{1}{c|}{1.03\%} & \multicolumn{1}{c|}{0.64\%} & \multicolumn{1}{c|}{2.21\%} & \multicolumn{1}{c|}{} & \multicolumn{1}{c|}{1E-04} & \multicolumn{1}{c|}{1.91\%} & \multicolumn{1}{c|}{1.24\%} & \multicolumn{1}{c|}{5.90\%} \\ \cline{1-4} \cline{6-9} \cline{11-14} 
\multicolumn{1}{|c|}{1E-05} & \multicolumn{1}{c|}{0.20\%} & \multicolumn{1}{c|}{0.15\%} & \multicolumn{1}{c|}{1.31\%} & \multicolumn{1}{c|}{} & \multicolumn{1}{c|}{1E-05} & \multicolumn{1}{c|}{0.40\%} & \multicolumn{1}{c|}{0.24\%} & \multicolumn{1}{c|}{0.82\%} & \multicolumn{1}{c|}{} & \multicolumn{1}{c|}{1E-05} & \multicolumn{1}{c|}{1.02\%} & \multicolumn{1}{c|}{0.66\%} & \multicolumn{1}{c|}{1.11\%} \\ \cline{1-4} \cline{6-9} \cline{11-14} 
\multicolumn{1}{|c|}{1E-06} & \multicolumn{1}{c|}{0.10\%} & \multicolumn{1}{c|}{0.07\%} & \multicolumn{1}{c|}{0.00\%} & \multicolumn{1}{c|}{} & \multicolumn{1}{c|}{1E-06} & \multicolumn{1}{c|}{0.26\%} & \multicolumn{1}{c|}{0.18\%} & \multicolumn{1}{c|}{0.43\%} & \multicolumn{1}{c|}{} & \multicolumn{1}{c|}{1E-06} & \multicolumn{1}{c|}{0.61\%} & \multicolumn{1}{c|}{0.49\%} & \multicolumn{1}{c|}{0.79\%} \\ \cline{1-4} \cline{6-9} \cline{11-14} 
\multicolumn{1}{|c|}{1E-07} & \multicolumn{1}{c|}{0.05\%} & \multicolumn{1}{c|}{0.04\%} & \multicolumn{1}{c|}{0.00\%} & \multicolumn{1}{c|}{} & \multicolumn{1}{c|}{1E-07} & \multicolumn{1}{c|}{0.10\%} & \multicolumn{1}{c|}{0.09\%} & \multicolumn{1}{c|}{0.12\%} & \multicolumn{1}{c|}{} & \multicolumn{1}{c|}{1E-07} & \multicolumn{1}{c|}{0.20\%} & \multicolumn{1}{c|}{0.20\%} & \multicolumn{1}{c|}{0.20\%} \\ \cline{1-4} \cline{6-9} \cline{11-14} 
\multicolumn{1}{|c|}{1E-08} & \multicolumn{1}{c|}{0.04\%} & \multicolumn{1}{c|}{0.03\%} & \multicolumn{1}{c|}{0.00\%} & \multicolumn{1}{c|}{} & \multicolumn{1}{c|}{1E-08} & \multicolumn{1}{c|}{0.05\%} & \multicolumn{1}{c|}{0.02\%} & \multicolumn{1}{c|}{0.06\%} & \multicolumn{1}{c|}{} & \multicolumn{1}{c|}{1E-08} & \multicolumn{1}{c|}{0.05\%} & \multicolumn{1}{c|}{0.04\%} & \multicolumn{1}{c|}{0.08\%} \\ \cline{1-4} \cline{6-9} \cline{11-14} 
 &  &  &  &  &  &  &  &  &  &  &  &  &  \\ \cline{1-4} \cline{6-9} \cline{11-14} 
\multicolumn{4}{|c|}{\textbf{Eastern Massachusetts}} & \multicolumn{1}{c|}{} & \multicolumn{4}{c|}{\textbf{Barcelona}} & \multicolumn{1}{c|}{} & \multicolumn{4}{c|}{\textbf{Berlin Center}} \\ \cline{1-4} \cline{6-9} \cline{11-14} 
\multicolumn{1}{|c|}{\textbf{Gap Level}} & \multicolumn{1}{c|}{\textbf{$\Delta$TSTT}} & \multicolumn{1}{c|}{\textbf{$\Delta$VMT}} & \multicolumn{1}{c|}{\textbf{PUL}} & \multicolumn{1}{c|}{} & \multicolumn{1}{c|}{\textbf{Gap Level}} & \multicolumn{1}{c|}{\textbf{$\Delta$TSTT}} & \multicolumn{1}{c|}{\textbf{$\Delta$VMT}} & \multicolumn{1}{c|}{\textbf{PUL}} & \multicolumn{1}{c|}{} & \multicolumn{1}{c|}{\textbf{Gap Level}} & \multicolumn{1}{c|}{\textbf{$\Delta$TSTT}} & \multicolumn{1}{c|}{\textbf{$\Delta$VMT}} & \multicolumn{1}{c|}{\textbf{PUL}} \\ \cline{1-4} \cline{6-9} \cline{11-14} 
\multicolumn{1}{|c|}{1E-03} & \multicolumn{1}{c|}{1.73\%} & \multicolumn{1}{c|}{1.34\%} & \multicolumn{1}{c|}{5.03\%} & \multicolumn{1}{c|}{} & \multicolumn{1}{c|}{1E-03} & \multicolumn{1}{c|}{2.02\%} & \multicolumn{1}{c|}{1.74\%} & \multicolumn{1}{c|}{6.10\%} & \multicolumn{1}{c|}{} & \multicolumn{1}{c|}{1E-03} & \multicolumn{1}{c|}{2.97\%} & \multicolumn{1}{c|}{2.59\%} & \multicolumn{1}{c|}{12.29\%} \\ \cline{1-4} \cline{6-9} \cline{11-14} 
\multicolumn{1}{|c|}{1E-04} & \multicolumn{1}{c|}{0.75\%} & \multicolumn{1}{c|}{0.44\%} & \multicolumn{1}{c|}{2.82\%} & \multicolumn{1}{c|}{} & \multicolumn{1}{c|}{1E-04} & \multicolumn{1}{c|}{1.00\%} & \multicolumn{1}{c|}{1.03\%} & \multicolumn{1}{c|}{2.88\%} & \multicolumn{1}{c|}{} & \multicolumn{1}{c|}{1E-04} & \multicolumn{1}{c|}{1.98\%} & \multicolumn{1}{c|}{1.32\%} & \multicolumn{1}{c|}{5.44\%} \\ \cline{1-4} \cline{6-9} \cline{11-14} 
\multicolumn{1}{|c|}{1E-05} & \multicolumn{1}{c|}{0.31\%} & \multicolumn{1}{c|}{0.10\%} & \multicolumn{1}{c|}{0.92\%} & \multicolumn{1}{c|}{} & \multicolumn{1}{c|}{1E-05} & \multicolumn{1}{c|}{0.52\%} & \multicolumn{1}{c|}{0.46\%} & \multicolumn{1}{c|}{1.02\%} & \multicolumn{1}{c|}{} & \multicolumn{1}{c|}{1E-05} & \multicolumn{1}{c|}{1.13\%} & \multicolumn{1}{c|}{0.74\%} & \multicolumn{1}{c|}{1.41\%} \\ \cline{1-4} \cline{6-9} \cline{11-14} 
\multicolumn{1}{|c|}{1E-06} & \multicolumn{1}{c|}{0.21\%} & \multicolumn{1}{c|}{0.07\%} & \multicolumn{1}{c|}{0.70\%} & \multicolumn{1}{c|}{} & \multicolumn{1}{c|}{1E-06} & \multicolumn{1}{c|}{0.24\%} & \multicolumn{1}{c|}{0.36\%} & \multicolumn{1}{c|}{0.73\%} & \multicolumn{1}{c|}{} & \multicolumn{1}{c|}{1E-06} & \multicolumn{1}{c|}{0.77\%} & \multicolumn{1}{c|}{0.52\%} & \multicolumn{1}{c|}{0.91\%} \\ \cline{1-4} \cline{6-9} \cline{11-14} 
\multicolumn{1}{|c|}{1E-07} & \multicolumn{1}{c|}{0.05\%} & \multicolumn{1}{c|}{0.03\%} & \multicolumn{1}{c|}{0.40\%} & \multicolumn{1}{c|}{} & \multicolumn{1}{c|}{1E-07} & \multicolumn{1}{c|}{0.13\%} & \multicolumn{1}{c|}{0.11\%} & \multicolumn{1}{c|}{0.32\%} & \multicolumn{1}{c|}{} & \multicolumn{1}{c|}{1E-07} & \multicolumn{1}{c|}{0.22\%} & \multicolumn{1}{c|}{0.18\%} & \multicolumn{1}{c|}{0.40\%} \\ \cline{1-4} \cline{6-9} \cline{11-14} 
\multicolumn{1}{|c|}{1E-08} & \multicolumn{1}{c|}{0.03\%} & \multicolumn{1}{c|}{0.01\%} & \multicolumn{1}{c|}{0.04\%} & \multicolumn{1}{c|}{} & \multicolumn{1}{c|}{1E-08} & \multicolumn{1}{c|}{0.02\%} & \multicolumn{1}{c|}{0.03\%} & \multicolumn{1}{c|}{0.08\%} & \multicolumn{1}{c|}{} & \multicolumn{1}{c|}{1E-08} & \multicolumn{1}{c|}{0.09\%} & \multicolumn{1}{c|}{0.04\%} & \multicolumn{1}{c|}{0.10\%} \\ \cline{1-4} \cline{6-9} \cline{11-14} 
 &  &  &  &  &  &  &  &  &  &  &  &  &  \\ \cline{1-4} \cline{6-9} \cline{11-14} 
\multicolumn{4}{|c|}{\textbf{Anaheim}} & \multicolumn{1}{c|}{} & \multicolumn{4}{c|}{\textbf{Winnipeg}} & \multicolumn{1}{c|}{} & \multicolumn{4}{c|}{\textbf{Chicago-Regional}} \\ \cline{1-4} \cline{6-9} \cline{11-14} 
\multicolumn{1}{|c|}{\textbf{Gap Level}} & \multicolumn{1}{c|}{\textbf{$\Delta$TSTT}} & \multicolumn{1}{c|}{\textbf{$\Delta$VMT}} & \multicolumn{1}{c|}{\textbf{PUL}} & \multicolumn{1}{c|}{} & \multicolumn{1}{c|}{\textbf{Gap Level}} & \multicolumn{1}{c|}{\textbf{$\Delta$TSTT}} & \multicolumn{1}{c|}{\textbf{$\Delta$VMT}} & \multicolumn{1}{c|}{\textbf{PUL}} & \multicolumn{1}{c|}{} & \multicolumn{1}{c|}{\textbf{Gap Level}} & \multicolumn{1}{c|}{\textbf{$\Delta$TSTT}} & \multicolumn{1}{c|}{\textbf{$\Delta$VMT}} & \multicolumn{1}{c|}{\textbf{PUL}} \\ \cline{1-4} \cline{6-9} \cline{11-14} 
\multicolumn{1}{|c|}{1E-03} & \multicolumn{1}{c|}{2.56\%} & \multicolumn{1}{c|}{1.45\%} & \multicolumn{1}{c|}{4.93\%} & \multicolumn{1}{c|}{} & \multicolumn{1}{c|}{1E-03} & \multicolumn{1}{c|}{2.97\%} & \multicolumn{1}{c|}{1.99\%} & \multicolumn{1}{c|}{6.02\%} & \multicolumn{1}{c|}{} & \multicolumn{1}{c|}{1E-03} & \multicolumn{1}{c|}{2.61\%} & \multicolumn{1}{c|}{2.38\%} & \multicolumn{1}{c|}{13.07\%} \\ \cline{1-4} \cline{6-9} \cline{11-14} 
\multicolumn{1}{|c|}{1E-04} & \multicolumn{1}{c|}{0.80\%} & \multicolumn{1}{c|}{0.60\%} & \multicolumn{1}{c|}{2.55\%} & \multicolumn{1}{c|}{} & \multicolumn{1}{c|}{1E-04} & \multicolumn{1}{c|}{1.01\%} & \multicolumn{1}{c|}{1.00\%} & \multicolumn{1}{c|}{2.45\%} & \multicolumn{1}{c|}{} & \multicolumn{1}{c|}{1E-04} & \multicolumn{1}{c|}{1.59\%} & \multicolumn{1}{c|}{1.10\%} & \multicolumn{1}{c|}{5.97\%} \\ \cline{1-4} \cline{6-9} \cline{11-14} 
\multicolumn{1}{|c|}{1E-05} & \multicolumn{1}{c|}{0.40\%} & \multicolumn{1}{c|}{0.30\%} & \multicolumn{1}{c|}{0.77\%} & \multicolumn{1}{c|}{} & \multicolumn{1}{c|}{1E-05} & \multicolumn{1}{c|}{0.53\%} & \multicolumn{1}{c|}{0.35\%} & \multicolumn{1}{c|}{0.95\%} & \multicolumn{1}{c|}{} & \multicolumn{1}{c|}{1E-05} & \multicolumn{1}{c|}{0.92\%} & \multicolumn{1}{c|}{0.61\%} & \multicolumn{1}{c|}{1.00\%} \\ \cline{1-4} \cline{6-9} \cline{11-14} 
\multicolumn{1}{|c|}{1E-06} & \multicolumn{1}{c|}{0.26\%} & \multicolumn{1}{c|}{0.19\%} & \multicolumn{1}{c|}{0.23\%} & \multicolumn{1}{c|}{} & \multicolumn{1}{c|}{1E-06} & \multicolumn{1}{c|}{0.29\%} & \multicolumn{1}{c|}{0.33\%} & \multicolumn{1}{c|}{0.52\%} & \multicolumn{1}{c|}{} & \multicolumn{1}{c|}{1E-06} & \multicolumn{1}{c|}{0.76\%} & \multicolumn{1}{c|}{0.54\%} & \multicolumn{1}{c|}{0.79\%} \\ \cline{1-4} \cline{6-9} \cline{11-14} 
\multicolumn{1}{|c|}{1E-07} & \multicolumn{1}{c|}{0.10\%} & \multicolumn{1}{c|}{0.08\%} & \multicolumn{1}{c|}{0.04\%} & \multicolumn{1}{c|}{} & \multicolumn{1}{c|}{1E-07} & \multicolumn{1}{c|}{0.12\%} & \multicolumn{1}{c|}{0.11\%} & \multicolumn{1}{c|}{0.18\%} & \multicolumn{1}{c|}{} & \multicolumn{1}{c|}{1E-07} & \multicolumn{1}{c|}{0.50\%} & \multicolumn{1}{c|}{0.23\%} & \multicolumn{1}{c|}{0.00\%} \\ \cline{1-4} \cline{6-9} \cline{11-14} 
\multicolumn{1}{|c|}{1E-08} & \multicolumn{1}{c|}{0.05\%} & \multicolumn{1}{c|}{0.01\%} & \multicolumn{1}{c|}{0.02\%} & \multicolumn{1}{c|}{} & \multicolumn{1}{c|}{1E-08} & \multicolumn{1}{c|}{0.02\%} & \multicolumn{1}{c|}{0.03\%} & \multicolumn{1}{c|}{0.05\%} & \multicolumn{1}{c|}{} & \multicolumn{1}{c|}{1E-08} & \multicolumn{1}{c|}{0.15\%} & \multicolumn{1}{c|}{0.10\%} & \multicolumn{1}{c|}{0.00\%} \\ \cline{1-4} \cline{6-9} \cline{11-14} 
 &  &  &  &  &  &  &  &  &  &  &  &  &  \\ \cline{1-4} \cline{6-9} \cline{11-14} 
\multicolumn{4}{|c|}{\textbf{Chicago Sketch}} & \multicolumn{1}{c|}{} & \multicolumn{4}{c|}{\textbf{Terrassa}} & \multicolumn{1}{c|}{} & \multicolumn{4}{c|}{\textbf{Philadelphia}} \\ \cline{1-4} \cline{6-9} \cline{11-14} 
\multicolumn{1}{|c|}{\textbf{Gap Level}} & \multicolumn{1}{c|}{\textbf{$\Delta$TSTT}} & \multicolumn{1}{c|}{\textbf{$\Delta$VMT}} & \multicolumn{1}{c|}{\textbf{PUL}} & \multicolumn{1}{c|}{} & \multicolumn{1}{c|}{\textbf{Gap Level}} & \multicolumn{1}{c|}{\textbf{$\Delta$TSTT}} & \multicolumn{1}{c|}{\textbf{$\Delta$VMT}} & \multicolumn{1}{c|}{\textbf{PUL}} & \multicolumn{1}{c|}{} & \multicolumn{1}{c|}{\textbf{Gap Level}} & \multicolumn{1}{c|}{\textbf{$\Delta$TSTT}} & \multicolumn{1}{c|}{\textbf{$\Delta$VMT}} & \multicolumn{1}{c|}{\textbf{PUL}} \\ \cline{1-4} \cline{6-9} \cline{11-14} 
\multicolumn{1}{|c|}{1E-03} & \multicolumn{1}{c|}{2.03\%} & \multicolumn{1}{c|}{1.02\%} & \multicolumn{1}{c|}{5.55\%} & \multicolumn{1}{c|}{} & \multicolumn{1}{c|}{1E-03} & \multicolumn{1}{c|}{2.74\%} & \multicolumn{1}{c|}{2.23\%} & \multicolumn{1}{c|}{7.07\%} & \multicolumn{1}{c|}{} & \multicolumn{1}{c|}{1E-03} & \multicolumn{1}{c|}{3.95\%} & \multicolumn{1}{c|}{2.80\%} & \multicolumn{1}{c|}{17.88\%} \\ \cline{1-4} \cline{6-9} \cline{11-14} 
\multicolumn{1}{|c|}{1E-04} & \multicolumn{1}{c|}{0.99\%} & \multicolumn{1}{c|}{0.50\%} & \multicolumn{1}{c|}{2.85\%} & \multicolumn{1}{c|}{} & \multicolumn{1}{c|}{1E-04} & \multicolumn{1}{c|}{1.15\%} & \multicolumn{1}{c|}{1.00\%} & \multicolumn{1}{c|}{3.07\%} & \multicolumn{1}{c|}{} & \multicolumn{1}{c|}{1E-04} & \multicolumn{1}{c|}{1.89\%} & \multicolumn{1}{c|}{1.50\%} & \multicolumn{1}{c|}{7.91\%} \\ \cline{1-4} \cline{6-9} \cline{11-14} 
\multicolumn{1}{|c|}{1E-05} & \multicolumn{1}{c|}{0.37\%} & \multicolumn{1}{c|}{0.25\%} & \multicolumn{1}{c|}{0.99\%} & \multicolumn{1}{c|}{} & \multicolumn{1}{c|}{1E-05} & \multicolumn{1}{c|}{0.52\%} & \multicolumn{1}{c|}{0.40\%} & \multicolumn{1}{c|}{1.01\%} & \multicolumn{1}{c|}{} & \multicolumn{1}{c|}{1E-05} & \multicolumn{1}{c|}{1.48\%} & \multicolumn{1}{c|}{1.26\%} & \multicolumn{1}{c|}{1.49\%} \\ \cline{1-4} \cline{6-9} \cline{11-14} 
\multicolumn{1}{|c|}{1E-06} & \multicolumn{1}{c|}{0.29\%} & \multicolumn{1}{c|}{0.09\%} & \multicolumn{1}{c|}{0.34\%} & \multicolumn{1}{c|}{} & \multicolumn{1}{c|}{1E-06} & \multicolumn{1}{c|}{0.29\%} & \multicolumn{1}{c|}{0.38\%} & \multicolumn{1}{c|}{0.51\%} & \multicolumn{1}{c|}{} & \multicolumn{1}{c|}{1E-06} & \multicolumn{1}{c|}{1.01\%} & \multicolumn{1}{c|}{0.76\%} & \multicolumn{1}{c|}{0.89\%} \\ \cline{1-4} \cline{6-9} \cline{11-14} 
\multicolumn{1}{|c|}{1E-07} & \multicolumn{1}{c|}{0.08\%} & \multicolumn{1}{c|}{0.04\%} & \multicolumn{1}{c|}{0.20\%} & \multicolumn{1}{c|}{} & \multicolumn{1}{c|}{1E-07} & \multicolumn{1}{c|}{0.13\%} & \multicolumn{1}{c|}{0.09\%} & \multicolumn{1}{c|}{0.20\%} & \multicolumn{1}{c|}{} & \multicolumn{1}{c|}{1E-07} & \multicolumn{1}{c|}{0.64\%} & \multicolumn{1}{c|}{0.45\%} & \multicolumn{1}{c|}{0.49\%} \\ \cline{1-4} \cline{6-9} \cline{11-14} 
\multicolumn{1}{|c|}{1E-08} & \multicolumn{1}{c|}{0.04\%} & \multicolumn{1}{c|}{0.02\%} & \multicolumn{1}{c|}{0.08\%} & \multicolumn{1}{c|}{} & \multicolumn{1}{c|}{1E-08} & \multicolumn{1}{c|}{0.02\%} & \multicolumn{1}{c|}{0.03\%} & \multicolumn{1}{c|}{0.08\%} & \multicolumn{1}{c|}{} & \multicolumn{1}{c|}{1E-08} & \multicolumn{1}{c|}{0.30\%} & \multicolumn{1}{c|}{0.05\%} & \multicolumn{1}{c|}{0.20\%} \\ \cline{1-4} \cline{6-9} \cline{11-14} 
\end{tabular}
\end{table}

\pagebreak

\begin{table}[]
\centering
\caption{Metric stabilization behavior comparison between Algorithm B and Gradient Projection}
\label{tab:algB_GP}
\begin{tabular}{ccccccclccccccc}
\cline{1-7} \cline{9-15}
\multicolumn{7}{|c|}{\textbf{Sioux Falls}} & \multicolumn{1}{l|}{} & \multicolumn{7}{c|}{\textbf{Winnipeg}} \\ \cline{1-7} \cline{9-15} 
\multicolumn{1}{|c|}{\multirow{2}{*}{\textbf{Gap Level}}} & \multicolumn{2}{c|}{\textbf{$\Delta$TSTT}} & \multicolumn{2}{c|}{\textbf{$\Delta$VMT}} & \multicolumn{2}{c|}{\textbf{PUL}} & \multicolumn{1}{l|}{} & \multicolumn{1}{c|}{\multirow{2}{*}{\textbf{Gap Level}}} & \multicolumn{2}{c|}{\textbf{$\Delta$TSTT}} & \multicolumn{2}{c|}{\textbf{$\Delta$VMT}} & \multicolumn{2}{c|}{\textbf{PUL}} \\ \cline{2-7} \cline{10-15} 
\multicolumn{1}{|c|}{} & \multicolumn{1}{c|}{\textbf{Alg-B}} & \multicolumn{1}{c|}{\textbf{GP}} & \multicolumn{1}{c|}{\textbf{Alg-B}} & \multicolumn{1}{c|}{\textbf{GP}} & \multicolumn{1}{c|}{\textbf{Alg-B}} & \multicolumn{1}{c|}{\textbf{GP}} & \multicolumn{1}{l|}{} & \multicolumn{1}{c|}{} & \multicolumn{1}{c|}{\textbf{Alg-B}} & \multicolumn{1}{c|}{\textbf{GP}} & \multicolumn{1}{c|}{\textbf{Alg-B}} & \multicolumn{1}{c|}{\textbf{GP}} & \multicolumn{1}{c|}{\textbf{Alg-B}} & \multicolumn{1}{c|}{\textbf{GP}} \\ \cline{1-7} \cline{9-15} 
\multicolumn{1}{|c|}{1E-03} & \multicolumn{1}{c|}{2.00\%} & \multicolumn{1}{c|}{2.13\%} & \multicolumn{1}{c|}{1.44\%} & \multicolumn{1}{c|}{1.48\%} & \multicolumn{1}{c|}{5.26\%} & \multicolumn{1}{c|}{5.26\%} & \multicolumn{1}{l|}{} & \multicolumn{1}{c|}{1E-03} & \multicolumn{1}{c|}{2.97\%} & \multicolumn{1}{c|}{2.84\%} & \multicolumn{1}{c|}{1.99\%} & \multicolumn{1}{c|}{1.81\%} & \multicolumn{1}{c|}{6.02\%} & \multicolumn{1}{c|}{5.49\%} \\ \cline{1-7} \cline{9-15} 
\multicolumn{1}{|c|}{1E-04} & \multicolumn{1}{c|}{0.80\%} & \multicolumn{1}{c|}{0.89\%} & \multicolumn{1}{c|}{0.59\%} & \multicolumn{1}{c|}{0.65\%} & \multicolumn{1}{c|}{2.63\%} & \multicolumn{1}{c|}{2.63\%} & \multicolumn{1}{l|}{} & \multicolumn{1}{c|}{1E-04} & \multicolumn{1}{c|}{1.01\%} & \multicolumn{1}{c|}{0.94\%} & \multicolumn{1}{c|}{1.00\%} & \multicolumn{1}{c|}{1.01\%} & \multicolumn{1}{c|}{2.45\%} & \multicolumn{1}{c|}{2.20\%} \\ \cline{1-7} \cline{9-15} 
\multicolumn{1}{|c|}{1E-05} & \multicolumn{1}{c|}{0.20\%} & \multicolumn{1}{c|}{0.35\%} & \multicolumn{1}{c|}{0.15\%} & \multicolumn{1}{c|}{0.18\%} & \multicolumn{1}{c|}{1.31\%} & \multicolumn{1}{c|}{1.31\%} & \multicolumn{1}{l|}{} & \multicolumn{1}{c|}{1E-05} & \multicolumn{1}{c|}{0.53\%} & \multicolumn{1}{c|}{0.47\%} & \multicolumn{1}{c|}{0.35\%} & \multicolumn{1}{c|}{0.30\%} & \multicolumn{1}{c|}{0.95\%} & \multicolumn{1}{c|}{0.91\%} \\ \cline{1-7} \cline{9-15} 
\multicolumn{1}{|c|}{1E-06} & \multicolumn{1}{c|}{0.10\%} & \multicolumn{1}{c|}{0.07\%} & \multicolumn{1}{c|}{0.07\%} & \multicolumn{1}{c|}{0.03\%} & \multicolumn{1}{c|}{0.00\%} & \multicolumn{1}{c|}{0.00\%} & \multicolumn{1}{l|}{} & \multicolumn{1}{c|}{1E-06} & \multicolumn{1}{c|}{0.29\%} & \multicolumn{1}{c|}{0.26\%} & \multicolumn{1}{c|}{0.33\%} & \multicolumn{1}{c|}{0.30\%} & \multicolumn{1}{c|}{0.52\%} & \multicolumn{1}{c|}{0.52\%} \\ \cline{1-7} \cline{9-15} 
\multicolumn{1}{|c|}{1E-07} & \multicolumn{1}{c|}{0.05\%} & \multicolumn{1}{c|}{0.02\%} & \multicolumn{1}{c|}{0.04\%} & \multicolumn{1}{c|}{0.01\%} & \multicolumn{1}{c|}{0.00\%} & \multicolumn{1}{c|}{0.00\%} & \multicolumn{1}{l|}{} & \multicolumn{1}{c|}{1E-07} & \multicolumn{1}{c|}{0.12\%} & \multicolumn{1}{c|}{0.11\%} & \multicolumn{1}{c|}{0.11\%} & \multicolumn{1}{c|}{0.10\%} & \multicolumn{1}{c|}{0.18\%} & \multicolumn{1}{c|}{0.16\%} \\ \cline{1-7} \cline{9-15} 
\multicolumn{1}{|c|}{1E-08} & \multicolumn{1}{c|}{0.04\%} & \multicolumn{1}{c|}{0.00\%} & \multicolumn{1}{c|}{0.03\%} & \multicolumn{1}{c|}{0.00\%} & \multicolumn{1}{c|}{0.00\%} & \multicolumn{1}{c|}{0.00\%} & \multicolumn{1}{l|}{} & \multicolumn{1}{c|}{1E-08} & \multicolumn{1}{c|}{0.02\%} & \multicolumn{1}{c|}{0.02\%} & \multicolumn{1}{c|}{0.03\%} & \multicolumn{1}{c|}{0.03\%} & \multicolumn{1}{c|}{0.05\%} & \multicolumn{1}{c|}{0.05\%} \\ \cline{1-7} \cline{9-15} 
\multicolumn{1}{l}{} & \multicolumn{1}{l}{} & \multicolumn{1}{l}{} & \multicolumn{1}{l}{} & \multicolumn{1}{l}{} & \multicolumn{1}{l}{} & \multicolumn{1}{l}{} &  & \multicolumn{1}{l}{} & \multicolumn{1}{l}{} & \multicolumn{1}{l}{} & \multicolumn{1}{l}{} & \multicolumn{1}{l}{} & \multicolumn{1}{l}{} & \multicolumn{1}{l}{} \\ \cline{1-7} \cline{9-15} 
\multicolumn{7}{|c|}{\textbf{Anaheim}} & \multicolumn{1}{l|}{} & \multicolumn{7}{c|}{\textbf{Austin}} \\ \cline{1-7} \cline{9-15} 
\multicolumn{1}{|c|}{\multirow{2}{*}{\textbf{Gap Level}}} & \multicolumn{2}{c|}{\textbf{$\Delta$TSTT}} & \multicolumn{2}{c|}{\textbf{$\Delta$VMT}} & \multicolumn{2}{c|}{\textbf{PUL}} & \multicolumn{1}{l|}{} & \multicolumn{1}{c|}{\multirow{2}{*}{\textbf{Gap Level}}} & \multicolumn{2}{c|}{\textbf{$\Delta$TSTT}} & \multicolumn{2}{c|}{\textbf{$\Delta$VMT}} & \multicolumn{2}{c|}{\textbf{PUL}} \\ \cline{2-7} \cline{10-15} 
\multicolumn{1}{|c|}{} & \multicolumn{1}{c|}{\textbf{Alg-B}} & \multicolumn{1}{c|}{\textbf{GP}} & \multicolumn{1}{c|}{\textbf{Alg-B}} & \multicolumn{1}{c|}{\textbf{GP}} & \multicolumn{1}{c|}{\textbf{Alg-B}} & \multicolumn{1}{c|}{\textbf{GP}} & \multicolumn{1}{l|}{} & \multicolumn{1}{c|}{} & \multicolumn{1}{c|}{\textbf{Alg-B}} & \multicolumn{1}{c|}{\textbf{GP}} & \multicolumn{1}{c|}{\textbf{Alg-B}} & \multicolumn{1}{c|}{\textbf{GP}} & \multicolumn{1}{c|}{\textbf{Alg-B}} & \multicolumn{1}{c|}{\textbf{GP}} \\ \cline{1-7} \cline{9-15} 
\multicolumn{1}{|c|}{1E-03} & \multicolumn{1}{c|}{2.56\%} & \multicolumn{1}{c|}{2.47\%} & \multicolumn{1}{c|}{1.45\%} & \multicolumn{1}{c|}{1.41\%} & \multicolumn{1}{c|}{4.93\%} & \multicolumn{1}{c|}{4.77\%} & \multicolumn{1}{l|}{} & \multicolumn{1}{c|}{1E-03} & \multicolumn{1}{c|}{2.87\%} & \multicolumn{1}{c|}{2.39\%} & \multicolumn{1}{c|}{2.47\%} & \multicolumn{1}{c|}{2.18\%} & \multicolumn{1}{c|}{11.24\%} & \multicolumn{1}{c|}{9.55\%} \\ \cline{1-7} \cline{9-15} 
\multicolumn{1}{|c|}{1E-04} & \multicolumn{1}{c|}{0.80\%} & \multicolumn{1}{c|}{0.77\%} & \multicolumn{1}{c|}{0.60\%} & \multicolumn{1}{c|}{0.57\%} & \multicolumn{1}{c|}{2.55\%} & \multicolumn{1}{c|}{2.58\%} & \multicolumn{1}{l|}{} & \multicolumn{1}{c|}{1E-04} & \multicolumn{1}{c|}{1.91\%} & \multicolumn{1}{c|}{1.59\%} & \multicolumn{1}{c|}{1.24\%} & \multicolumn{1}{c|}{1.09\%} & \multicolumn{1}{c|}{5.90\%} & \multicolumn{1}{c|}{5.38\%} \\ \cline{1-7} \cline{9-15} 
\multicolumn{1}{|c|}{1E-05} & \multicolumn{1}{c|}{0.40\%} & \multicolumn{1}{c|}{0.41\%} & \multicolumn{1}{c|}{0.30\%} & \multicolumn{1}{c|}{0.31\%} & \multicolumn{1}{c|}{0.77\%} & \multicolumn{1}{c|}{0.72\%} & \multicolumn{1}{l|}{} & \multicolumn{1}{c|}{1E-05} & \multicolumn{1}{c|}{1.02\%} & \multicolumn{1}{c|}{0.85\%} & \multicolumn{1}{c|}{0.66\%} & \multicolumn{1}{c|}{0.59\%} & \multicolumn{1}{c|}{1.11\%} & \multicolumn{1}{c|}{0.90\%} \\ \cline{1-7} \cline{9-15} 
\multicolumn{1}{|c|}{1E-06} & \multicolumn{1}{c|}{0.26\%} & \multicolumn{1}{c|}{0.25\%} & \multicolumn{1}{c|}{0.19\%} & \multicolumn{1}{c|}{0.18\%} & \multicolumn{1}{c|}{0.23\%} & \multicolumn{1}{c|}{0.24\%} & \multicolumn{1}{l|}{} & \multicolumn{1}{c|}{1E-06} & \multicolumn{1}{c|}{0.61\%} & \multicolumn{1}{c|}{0.58\%} & \multicolumn{1}{c|}{0.49\%} & \multicolumn{1}{c|}{0.39\%} & \multicolumn{1}{c|}{0.79\%} & \multicolumn{1}{c|}{0.72\%} \\ \cline{1-7} \cline{9-15} 
\multicolumn{1}{|c|}{1E-07} & \multicolumn{1}{c|}{0.10\%} & \multicolumn{1}{c|}{0.10\%} & \multicolumn{1}{c|}{0.08\%} & \multicolumn{1}{c|}{0.08\%} & \multicolumn{1}{c|}{0.04\%} & \multicolumn{1}{c|}{0.04\%} & \multicolumn{1}{l|}{} & \multicolumn{1}{c|}{1E-07} & \multicolumn{1}{c|}{0.20\%} & \multicolumn{1}{c|}{0.19\%} & \multicolumn{1}{c|}{0.20\%} & \multicolumn{1}{c|}{0.19\%} & \multicolumn{1}{c|}{0.20\%} & \multicolumn{1}{c|}{0.19\%} \\ \cline{1-7} \cline{9-15} 
\multicolumn{1}{|c|}{1E-08} & \multicolumn{1}{c|}{0.05\%} & \multicolumn{1}{c|}{0.05\%} & \multicolumn{1}{c|}{0.01\%} & \multicolumn{1}{c|}{0.01\%} & \multicolumn{1}{c|}{0.02\%} & \multicolumn{1}{c|}{0.02\%} & \multicolumn{1}{l|}{} & \multicolumn{1}{c|}{1E-08} & \multicolumn{1}{c|}{0.05\%} & \multicolumn{1}{c|}{0.05\%} & \multicolumn{1}{c|}{0.04\%} & \multicolumn{1}{c|}{0.03\%} & \multicolumn{1}{c|}{0.08\%} & \multicolumn{1}{c|}{0.08\%} \\ \cline{1-7} \cline{9-15} 
\multicolumn{1}{l}{} & \multicolumn{1}{l}{} & \multicolumn{1}{l}{} & \multicolumn{1}{l}{} & \multicolumn{1}{l}{} & \multicolumn{1}{l}{} & \multicolumn{1}{l}{} &  & \multicolumn{1}{l}{} & \multicolumn{1}{l}{} & \multicolumn{1}{l}{} & \multicolumn{1}{l}{} & \multicolumn{1}{l}{} & \multicolumn{1}{l}{} & \multicolumn{1}{l}{} \\ \cline{1-7} \cline{9-15} 
\multicolumn{7}{|c|}{\textbf{Chicago-Sketch}} & \multicolumn{1}{l|}{} & \multicolumn{7}{c|}{\textbf{Philadelphia}} \\ \cline{1-7} \cline{9-15} 
\multicolumn{1}{|c|}{\multirow{2}{*}{\textbf{Gap Level}}} & \multicolumn{2}{c|}{\textbf{$\Delta$TSTT}} & \multicolumn{2}{c|}{\textbf{$\Delta$VMT}} & \multicolumn{2}{c|}{\textbf{PUL}} & \multicolumn{1}{l|}{} & \multicolumn{1}{c|}{\multirow{2}{*}{\textbf{Gap Level}}} & \multicolumn{2}{c|}{\textbf{$\Delta$TSTT}} & \multicolumn{2}{c|}{\textbf{$\Delta$VMT}} & \multicolumn{2}{c|}{\textbf{PUL}} \\ \cline{2-7} \cline{10-15} 
\multicolumn{1}{|c|}{} & \multicolumn{1}{c|}{\textbf{Alg-B}} & \multicolumn{1}{c|}{\textbf{GP}} & \multicolumn{1}{c|}{\textbf{Alg-B}} & \multicolumn{1}{c|}{\textbf{GP}} & \multicolumn{1}{c|}{\textbf{Alg-B}} & \multicolumn{1}{c|}{\textbf{GP}} & \multicolumn{1}{l|}{} & \multicolumn{1}{c|}{} & \multicolumn{1}{c|}{\textbf{Alg-B}} & \multicolumn{1}{c|}{\textbf{GP}} & \multicolumn{1}{c|}{\textbf{Alg-B}} & \multicolumn{1}{c|}{\textbf{GP}} & \multicolumn{1}{c|}{\textbf{Alg-B}} & \multicolumn{1}{c|}{\textbf{GP}} \\ \cline{1-7} \cline{9-15} 
\multicolumn{1}{|c|}{1E-03} & \multicolumn{1}{c|}{2.03\%} & \multicolumn{1}{c|}{1.88\%} & \multicolumn{1}{c|}{1.02\%} & \multicolumn{1}{c|}{0.94\%} & \multicolumn{1}{c|}{5.55\%} & \multicolumn{1}{c|}{4.87\%} & \multicolumn{1}{l|}{} & \multicolumn{1}{c|}{1E-03} & \multicolumn{1}{c|}{3.95\%} & \multicolumn{1}{c|}{3.83\%} & \multicolumn{1}{c|}{2.80\%} & \multicolumn{1}{c|}{2.43\%} & \multicolumn{1}{c|}{17.88\%} & \multicolumn{1}{c|}{15.82\%} \\ \cline{1-7} \cline{9-15} 
\multicolumn{1}{|c|}{1E-04} & \multicolumn{1}{c|}{0.99\%} & \multicolumn{1}{c|}{0.89\%} & \multicolumn{1}{c|}{0.50\%} & \multicolumn{1}{c|}{0.45\%} & \multicolumn{1}{c|}{2.85\%} & \multicolumn{1}{c|}{2.49\%} & \multicolumn{1}{l|}{} & \multicolumn{1}{c|}{1E-04} & \multicolumn{1}{c|}{1.89\%} & \multicolumn{1}{c|}{1.62\%} & \multicolumn{1}{c|}{1.50\%} & \multicolumn{1}{c|}{1.24\%} & \multicolumn{1}{c|}{7.91\%} & \multicolumn{1}{c|}{7.53\%} \\ \cline{1-7} \cline{9-15} 
\multicolumn{1}{|c|}{1E-05} & \multicolumn{1}{c|}{0.37\%} & \multicolumn{1}{c|}{0.38\%} & \multicolumn{1}{c|}{0.25\%} & \multicolumn{1}{c|}{0.25\%} & \multicolumn{1}{c|}{0.99\%} & \multicolumn{1}{c|}{0.96\%} & \multicolumn{1}{l|}{} & \multicolumn{1}{c|}{1E-05} & \multicolumn{1}{c|}{1.48\%} & \multicolumn{1}{c|}{1.05\%} & \multicolumn{1}{c|}{1.26\%} & \multicolumn{1}{c|}{0.79\%} & \multicolumn{1}{c|}{1.49\%} & \multicolumn{1}{c|}{1.26\%} \\ \cline{1-7} \cline{9-15} 
\multicolumn{1}{|c|}{1E-06} & \multicolumn{1}{c|}{0.29\%} & \multicolumn{1}{c|}{0.29\%} & \multicolumn{1}{c|}{0.09\%} & \multicolumn{1}{c|}{0.08\%} & \multicolumn{1}{c|}{0.34\%} & \multicolumn{1}{c|}{0.32\%} & \multicolumn{1}{l|}{} & \multicolumn{1}{c|}{1E-06} & \multicolumn{1}{c|}{1.01\%} & \multicolumn{1}{c|}{0.63\%} & \multicolumn{1}{c|}{0.76\%} & \multicolumn{1}{c|}{0.35\%} & \multicolumn{1}{c|}{0.89\%} & \multicolumn{1}{c|}{0.82\%} \\ \cline{1-7} \cline{9-15} 
\multicolumn{1}{|c|}{1E-07} & \multicolumn{1}{c|}{0.08\%} & \multicolumn{1}{c|}{0.08\%} & \multicolumn{1}{c|}{0.04\%} & \multicolumn{1}{c|}{0.04\%} & \multicolumn{1}{c|}{0.20\%} & \multicolumn{1}{c|}{0.18\%} & \multicolumn{1}{l|}{} & \multicolumn{1}{c|}{1E-07} & \multicolumn{1}{c|}{0.64\%} & \multicolumn{1}{c|}{0.36\%} & \multicolumn{1}{c|}{0.45\%} & \multicolumn{1}{c|}{0.14\%} & \multicolumn{1}{c|}{0.49\%} & \multicolumn{1}{c|}{0.37\%} \\ \cline{1-7} \cline{9-15} 
\multicolumn{1}{|c|}{1E-08} & \multicolumn{1}{c|}{0.04\%} & \multicolumn{1}{c|}{0.04\%} & \multicolumn{1}{c|}{0.02\%} & \multicolumn{1}{c|}{0.02\%} & \multicolumn{1}{c|}{0.08\%} & \multicolumn{1}{c|}{0.07\%} & \multicolumn{1}{l|}{} & \multicolumn{1}{c|}{1E-08} & \multicolumn{1}{c|}{0.30\%} & \multicolumn{1}{c|}{0.10\%} & \multicolumn{1}{c|}{0.05\%} & \multicolumn{1}{c|}{0.03\%} & \multicolumn{1}{c|}{0.20\%} & \multicolumn{1}{c|}{0.18\%} \\ \cline{1-7} \cline{9-15} 
\end{tabular}
\end{table}

\end{landscape}

\clearpage

\newcommand{\fakesection}[1]{%
  \par\refstepcounter{section}
  \sectionmark{#1}
  \addcontentsline{toc}{section}{\protect\numberline{\thesection}#1}
}

\fakesection{Tables and Figures}

\renewcommand\thetable{\arabic{table}}  
\setcounter{table}{0}   
\renewcommand\thefigure{\arabic{figure}}  
\setcounter{figure}{0}

\begin{landscape}
\begin{table}
\caption{Toy network costs for various scenarios}
\begin{tabular}{||c|c|c|c|c|c|l|}
\hline
 & \textbf{Separable} & \textbf{Symmetric-full} & \textbf{Symmetric-partial} & \textbf{Asymmetric-full} & \multicolumn{2}{c|}{\textbf{Asymmetric-partial}} \\ \hline
\textbf{Link 1} & $15+x_1$ & $15+0.5x_1+0.167(x_2+x_3+x_4)$ & $15+0.75x_1+0.25x_2$ & $15+0.5x_1+0.15x_2+0.167x_3+0.183x_4$ & \multicolumn{2}{c|}{$15+0.75x_1+0.25x_2$} \\ \hline
\textbf{Link 2} & $10+x_2$ & $10+0.5x_2+0.167(x_1+x_3+x_4)$ & $10+0.75x_2+0.25x_1$ & $10+0.5x_2+0.167x_1+0.183x_3+0.15x_4$ & \multicolumn{2}{c|}{$10+0.75x_2+0.3x_1$} \\ \hline
\textbf{Link 3} & $10+x_3$ & $10+0.5x_3+0.167(x_1+x_2+x_4)$ & $10+0.75x_3+0.25x_4$ & $10+0.5x_3+0.183x_1+0.167x_2+0.15x_4$ & \multicolumn{2}{c|}{$10+0.75x_3+0.3x_4$} \\ \hline
\textbf{Link 4} & $15+x_4$ & $15+0.5x_4+0.167(x_1+x_2+x_3)$ & $15+0.75x_4+0.25x_3$ & $15+0.5x_4+0.15x_1+0.183x_2+0.167x_3$ & \multicolumn{2}{c|}{$15+0.75x_4+0.25x_3$} \\ \hline
\end{tabular}
\label{table:toy_example_cases}
\end{table}

\end{landscape}

\begin{landscape}

\begin{table}
\caption{Toy network convergence behavior}
\label{table:toy_numbers}
\begin{tabular}{|c|c|c|c|c|c|}
\hline
 & \multicolumn{5}{c|}{\textbf{Relative gap}} \\ \hline
\textbf{Iteration} & \textbf{Separable} & \textbf{Symmetric-full} & \textbf{Symmetric-partial} & \textbf{Asymmetric-full} & \textbf{Asymmetric-partial} \\ \hline
1 & 6.0000 & 1.0000 & 4.5000 & 0.9048 & 4.5000 \\ \hline
2 & 1.6670 & 0.2444 & 0.3542 & 0.3150 & 0.7500 \\ \hline
3 & 1.1806 & 0.1531 & 0.3399 & 0.1537 & 0.3632 \\ \hline
4 & 0.4289 & 0.0733 & 0.1482 & 0.0828 & 0.2465 \\ \hline
5 & 0.6968 & 0.0384 & 0.0828 & 0.0396 & 0.1591 \\ \hline
\textbf{} & \multicolumn{5}{c|}{\textbf{Flows}} \\ \hline
\textbf{Iteration} & \textbf{Separable} & \textbf{Symmetric-full} & \textbf{Symmetric-partial} & \textbf{Asymmetric-full} & \textbf{Asymmetric-partial} \\ \hline
1 & {[}0.00,60.00,0.00,0.00{]} & {[}0.00,60.00,0.00,0.00{]} & {[}0.00,60.00,0.00,0.00{]} & {[}0.00,60.00,0.00,0.00{]} & {[}0.00,60.00,0.00,0.00{]} \\ \hline
2 & {[}0.00,30.00,30.00,0.00{]} & {[}0.00,40.00,20.00,0.00{]} & {[}0.00,30.00,30.00,0.00{]} & {[}0.00,41.00,19.00,0.00{]} & {[}0.00,40.00,20.00,0.00{]} \\ \hline
3 & {[}25.00,17.50,17.50,0.00{]} & {[}10.00,31.67,18.33,0.00{]} & {[}11.33,24.33,24.33,0.00{]} & {[}0.00,31.34,15.98,12.68{]} & {[}0.00,26.67,16.67,16.67{]} \\ \hline
4 & {[}12.50,11.25,11.25,25.00{]} & {[}6.67,26.11,17.22,10.00{]} & {[}5.67,18.48,20.37,15.48{]} & {[}9.78,25.73,15.19,9.30{]} & {[}15.55,21.11,13.33,10.00{]} \\ \hline
5 & {[}9.38,26.88,8.12,15.62{]} & {[}5.18,24.63,22.78,7.41{]} & {[}14.18,18.16,17.46,10.20{]} & {[}6.91,22.65,23.22,7.21{]} & {[}9.26,16.30,26.67,7.77{]} \\ \hline
\end{tabular}
\end{table}
\end{landscape}

\begin{table}
{
\caption{Description of networks used}\label{table:network_data}
\noindent\makebox[\textwidth]{%
\begin{tabular}{||l r r r r||} 
    \hline
    Network name & Zones & Links & Nodes & Trips\\ [0.5ex] 
    \hline\hline
    SiouxFalls & 24 & 76 & 24 & 360,600\\
    Eastern-Massachusetts & 74 & 258 & 74 & 65,576\\
    Chicago-sketch & 387 & 2950 & 933 & 1,260,907 \\
    Barcelona & 110 & 2522 & 1020 & 184,679\\
    Chicago-Regional & 1790 & 39018 & 12982 & 1,360,427 \\
    \hline
\end{tabular}}
}
\end{table}

\clearpage

\begin{table}
{
\caption{Description of networks used}\label{table:nets_STAP}
\begin{center}
\noindent\makebox[\textwidth]{%
\begin{tabular}{||l r r r r c||} 
    \hline
    Network name & Zones & Links & Nodes & Trips & Average flow-to-capacity ratio\\ [0.5ex] 
    \hline\hline
    SiouxFalls & 24 & 76 & 24 & 360,600 & 1.612 \\
    Eastern-Massachusetts & 74 & 258 & 74 & 65,576 & 0.163\\
    Anaheim & 38 & 914 & 416 & 104,694 & 0.297 \\
    \hline
    Chicago-sketch & 387 & 2950 & 933 & 1,260,907 & 0.257 \\
    Berlin-Prenzlauerberg-Center & 98 & 2184 & 975 & 23,648 & 0.121 \\
    Barcelona & 110 & 2522 & 1020 & 184,679 & 1.137 \\
    Winnipeg & 147 & 2836 & 1052 & 64,784 & 2.028\\
    Terrassa & 55 & 3264 & 1609 & 25,225,700 & 5.964\\
    \hline
    Austin & 7388 & 18961 & 7388 & 739,351 & 0.875 \\
    Berlin-Center & 865 & 28376 & 12981 & 168,222 & 0.092 \\
    Chicago-Regional & 1790 & 39018 & 12982 & 1,360,427 & 0.522 \\
    Philadelphia & 1525 & 40003 & 13389 & 18,503,872 & 0.949 \\
    \hline
\end{tabular}}
\end{center}
}
\end{table}

\clearpage

\begin{figure}
\centering
\subfigure[Example network]{%
\resizebox*{7cm}{!}{\label{fig:STAP-net_mult_opt}\includegraphics{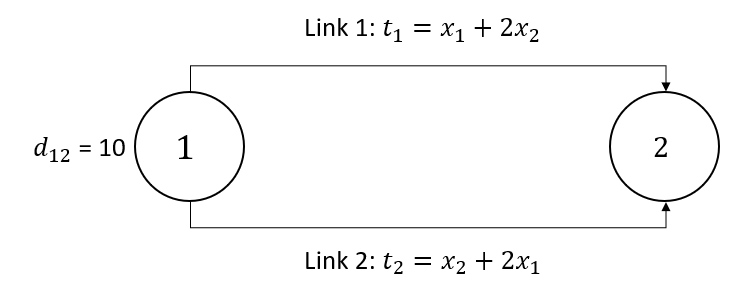}}}\hspace{5pt}
\subfigure[Objective function visualization]{%
\resizebox*{7cm}{!}{\label{fig:STAP_obj_fn}\includegraphics{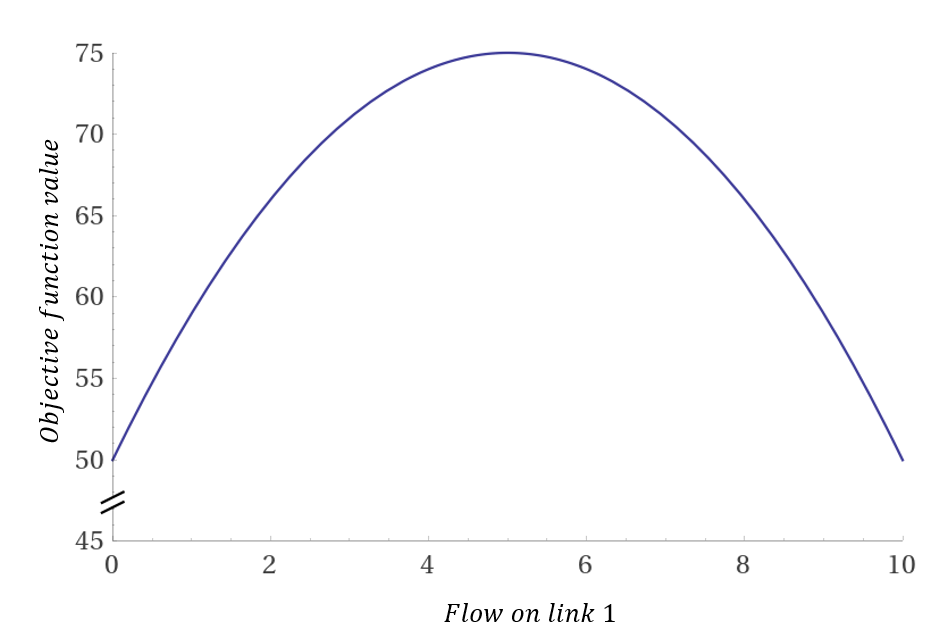}}}
\caption{S-TAP example with multiple extreme points} \label{fig:S-TAP_multiple_optima}
\end{figure}

\clearpage

\begin{figure}
\noindent\makebox[\textwidth]{%
\includegraphics[scale=0.6]{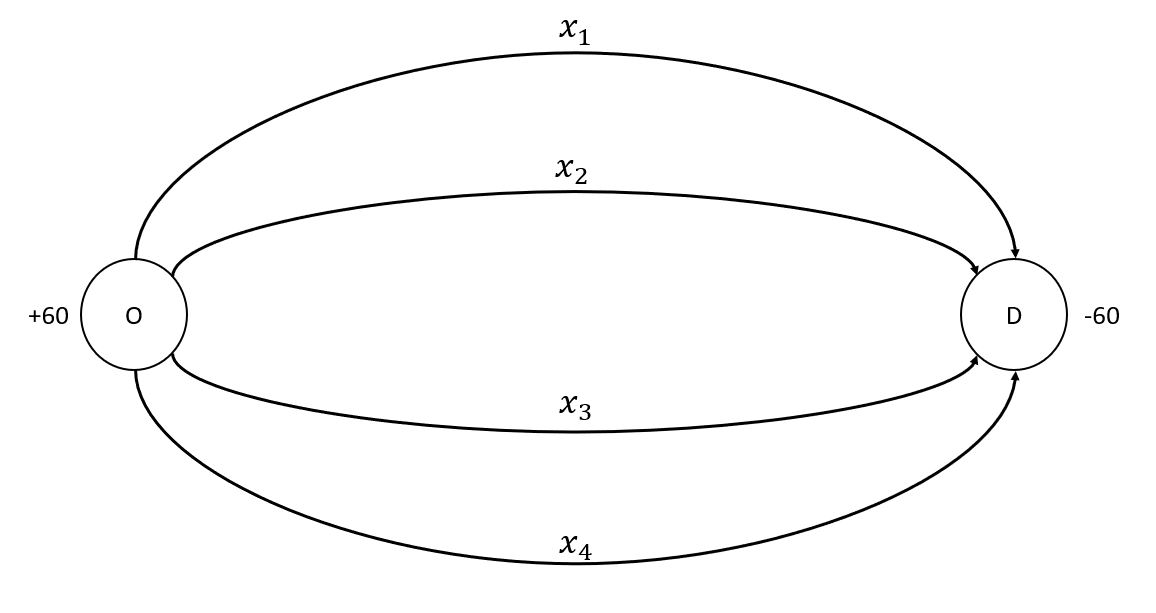}}
\centering
\caption{Toy network}
\label{fig:toy_example}
\centering
\end{figure}

\clearpage

\begin{figure}
    \centering
    \includegraphics[width=\textwidth]{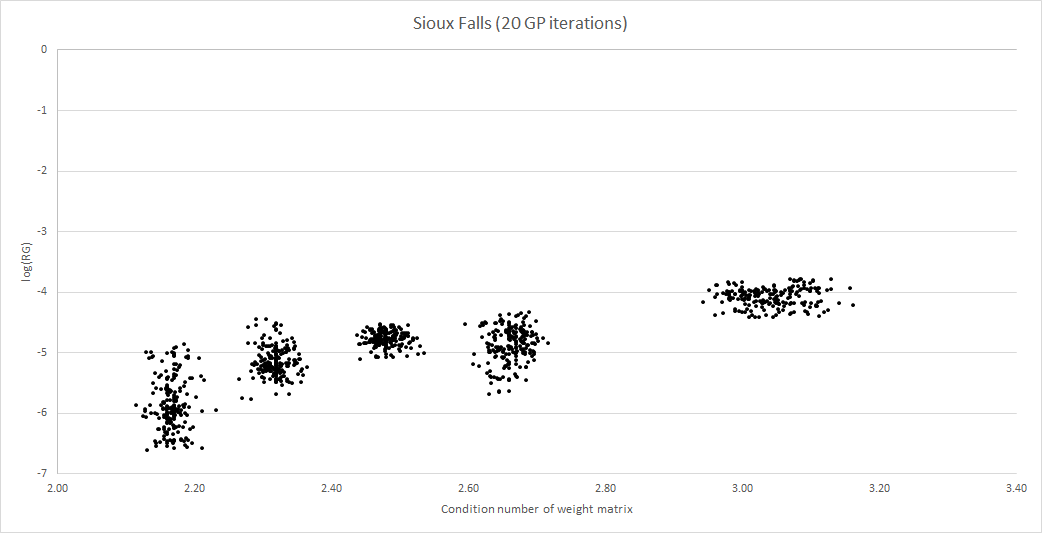}
    \caption{Condition number behavior of Sioux Falls problem instances}
    \label{fig:condition_num}
\end{figure}

\clearpage

\begin{figure}
\centering
\subfigure[Asymmetric to symmetric weight matrix convergence]{\label{fig:SF_a}\includegraphics[width=0.8\textwidth]{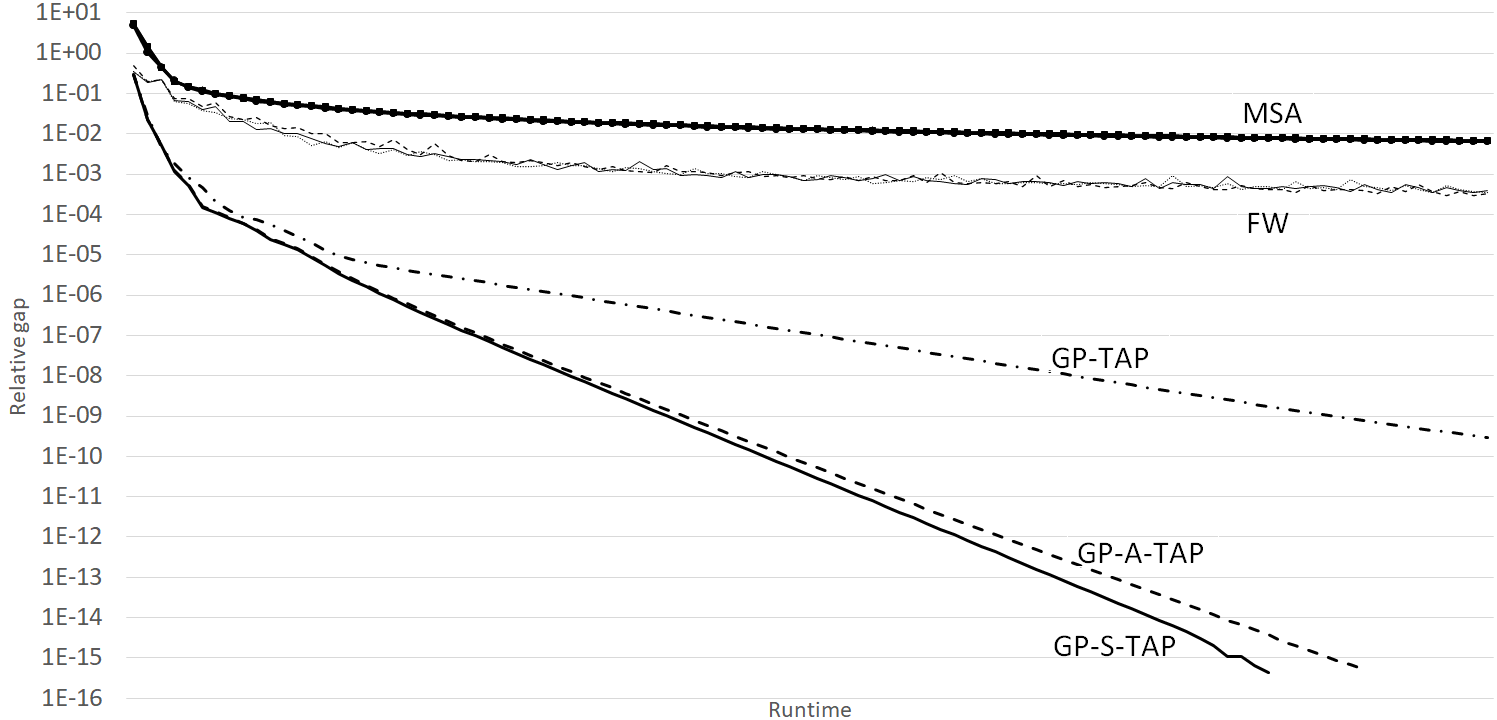}}
\subfigure[Convergence for S-TAP GP with $N$-link cost dependency]{\label{fig:SF_b}\includegraphics[width=0.8\textwidth]{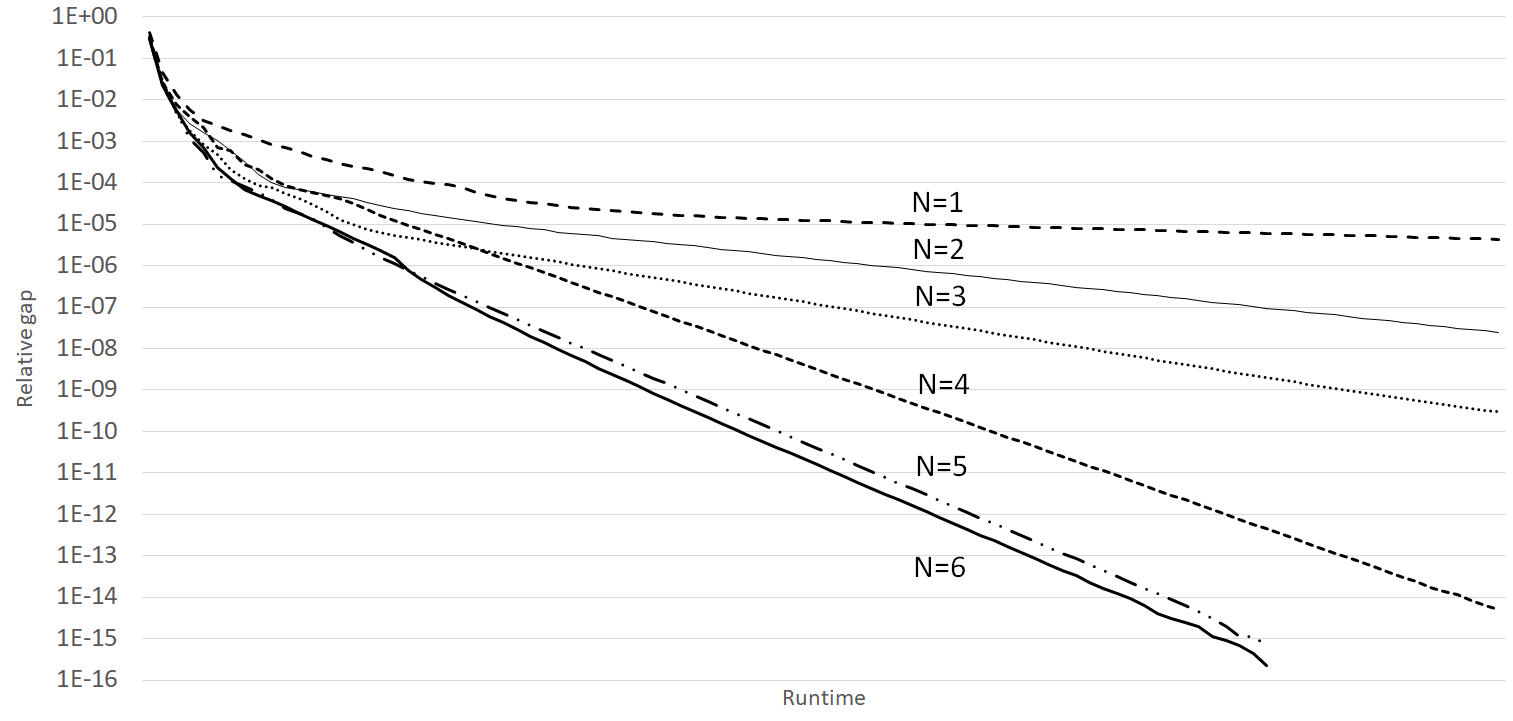}}
\subfigure[Convergence for A-TAP GP with $N$-link cost dependency]{\label{fig:SF_c}\includegraphics[width=0.8\textwidth]{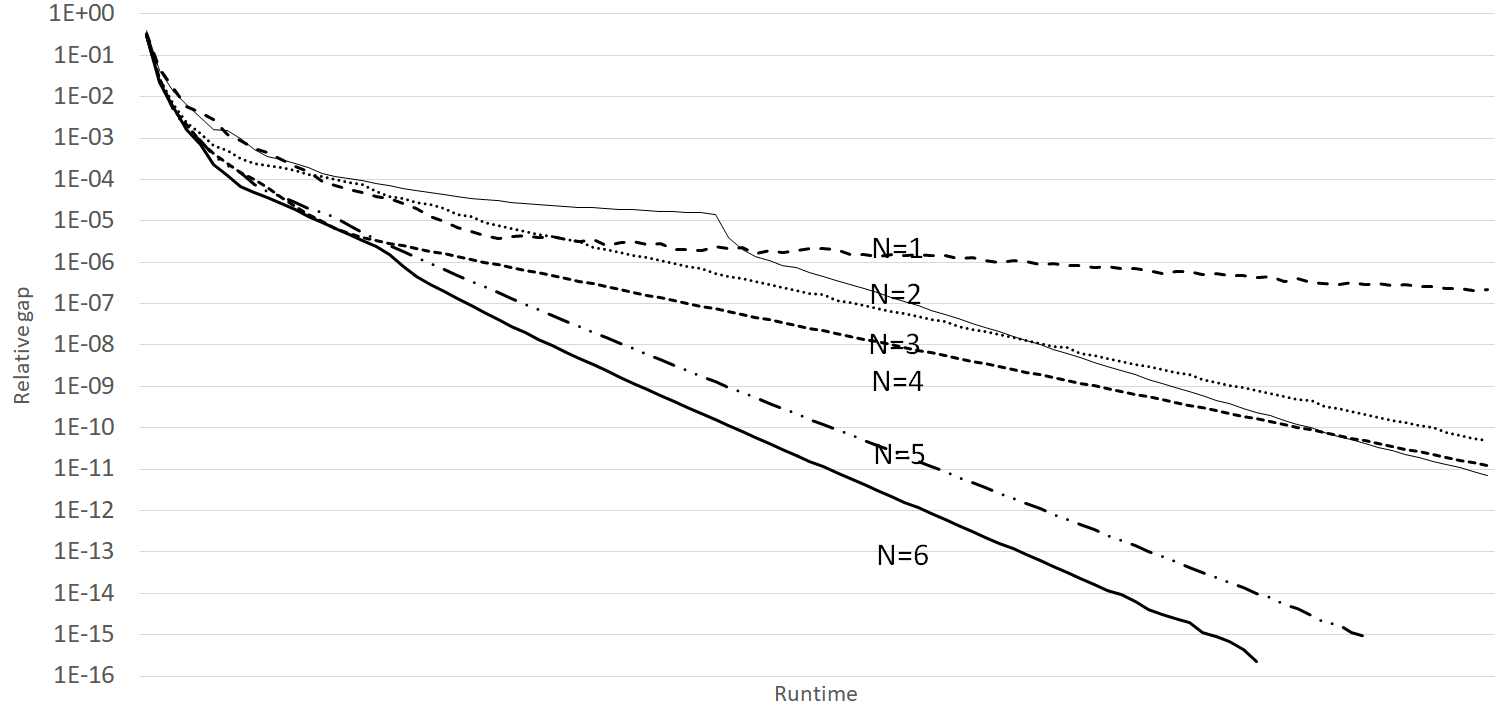}}
\caption{Experimental results for Sioux Falls network}
\label{fig:SF_compiled_results}
\end{figure}

\clearpage

\begin{figure}
    \centering
    \includegraphics[width=\textwidth]{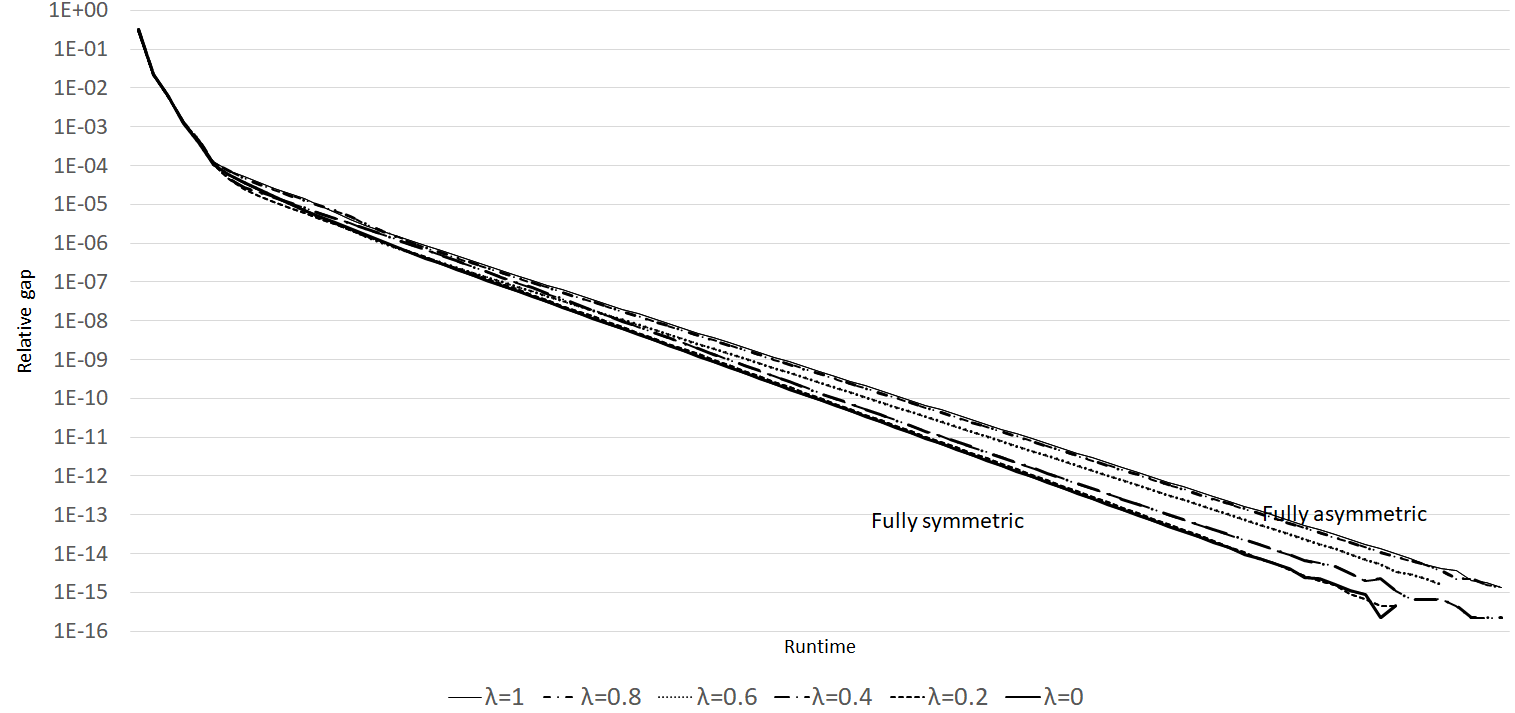}
    \caption{Sioux Falls asymmetric to symmetric weight matrix convergence}
    \label{fig:SF_asym_to_sym}
\end{figure}

\clearpage

\begin{figure}
\centering
\subfigure[Convergence for S-TAP GP with $N$-link cost dependency]{\label{fig:EMA_a}\includegraphics[width=0.8\textwidth]{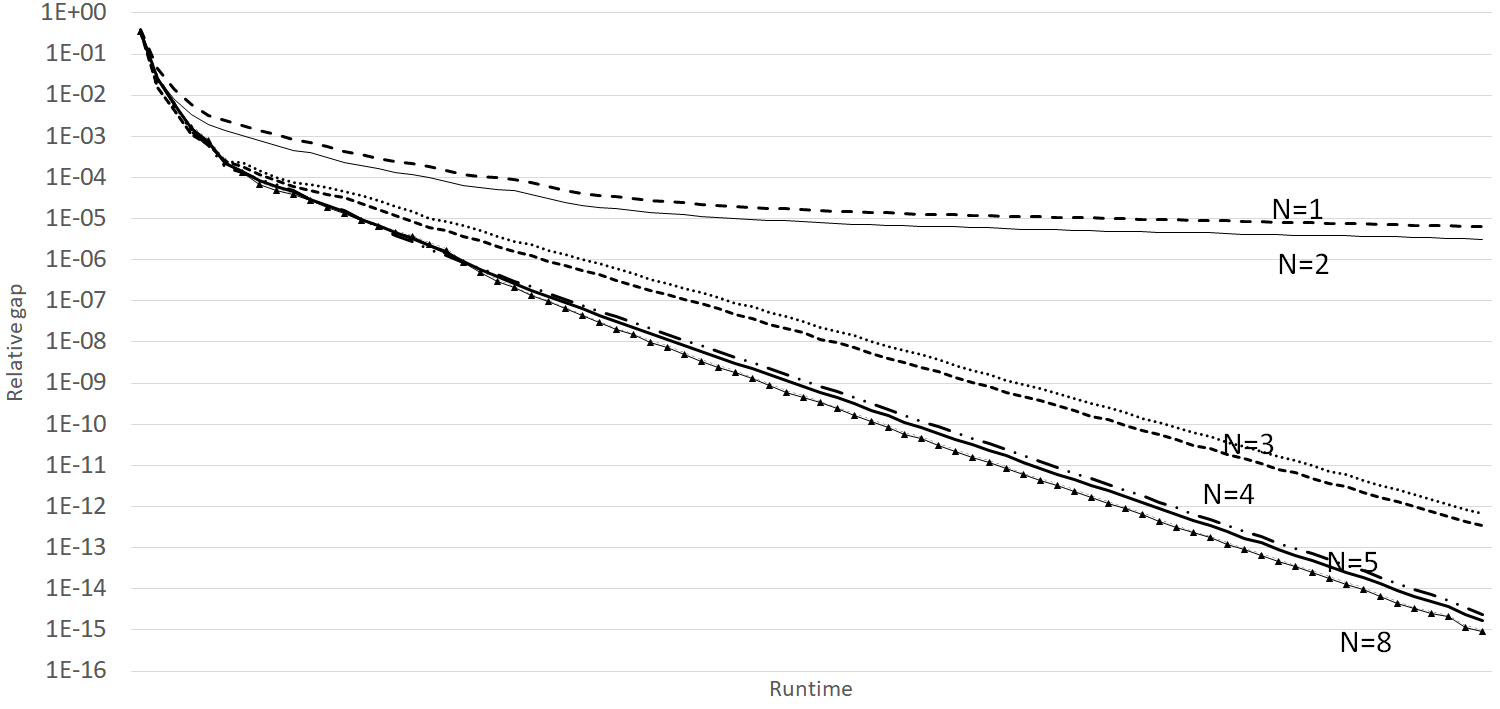}}
\subfigure[Convergence for A-TAP GP with $N$-link cost dependency]{\label{fig:EMA_b}\includegraphics[width=0.8\textwidth]{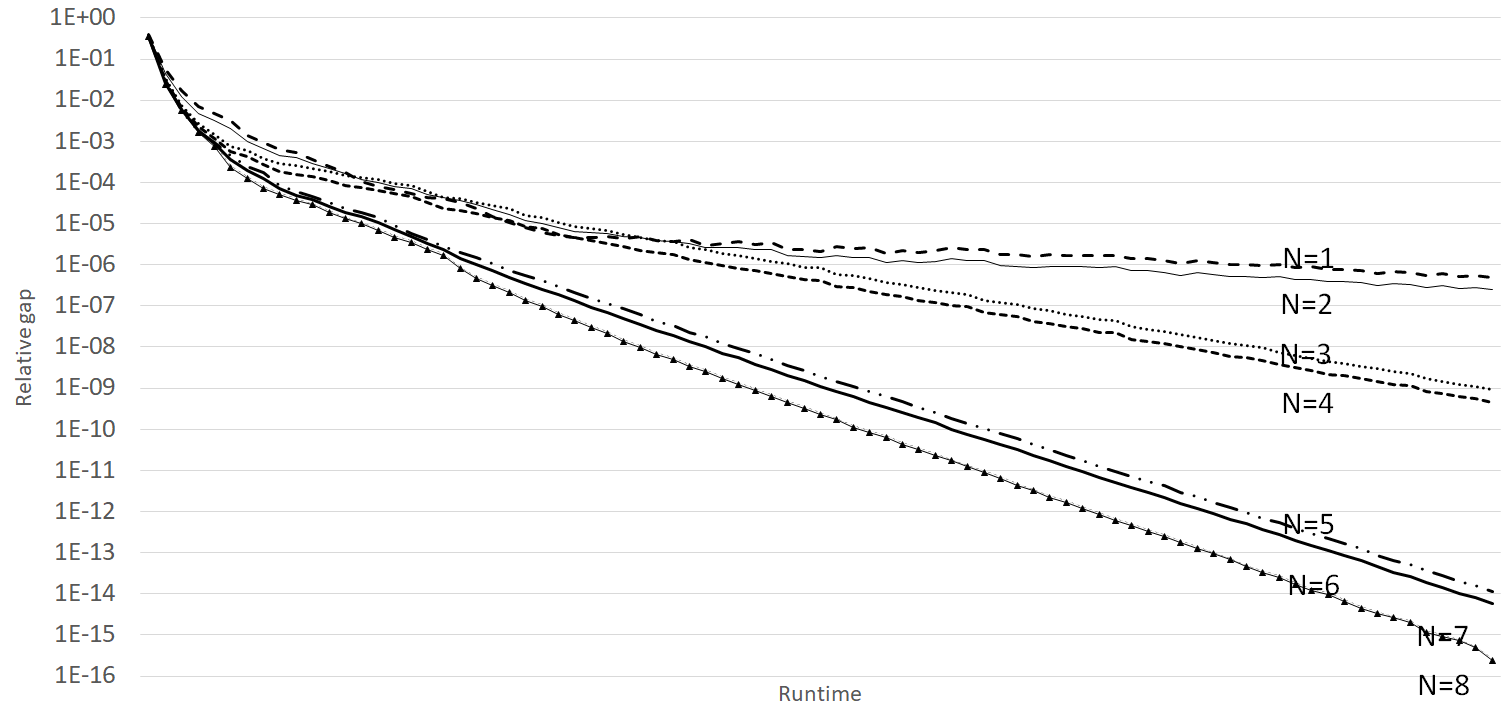}}
\subfigure[Asymmetric to symmetric weight matrix convergence]{\label{fig:EMA_c}\includegraphics[width=0.8\textwidth]{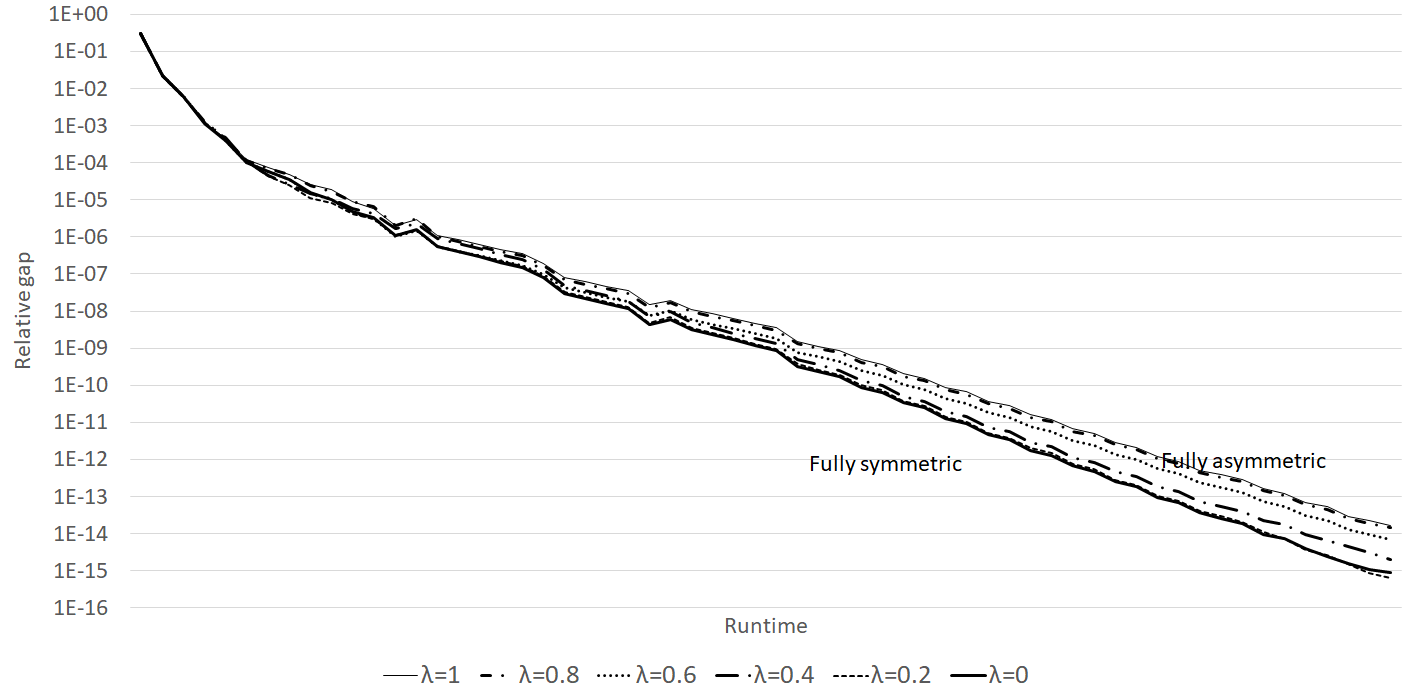}}
\caption{Experimental results for Eastern Massachusetts network}
\label{fig:EMA_compiled_results}
\end{figure}

\clearpage

\begin{figure}
\centering
\subfigure[Convergence for S-TAP GP with $N$-link cost dependency]{\label{fig:CS_a}\includegraphics[width=0.8\textwidth]{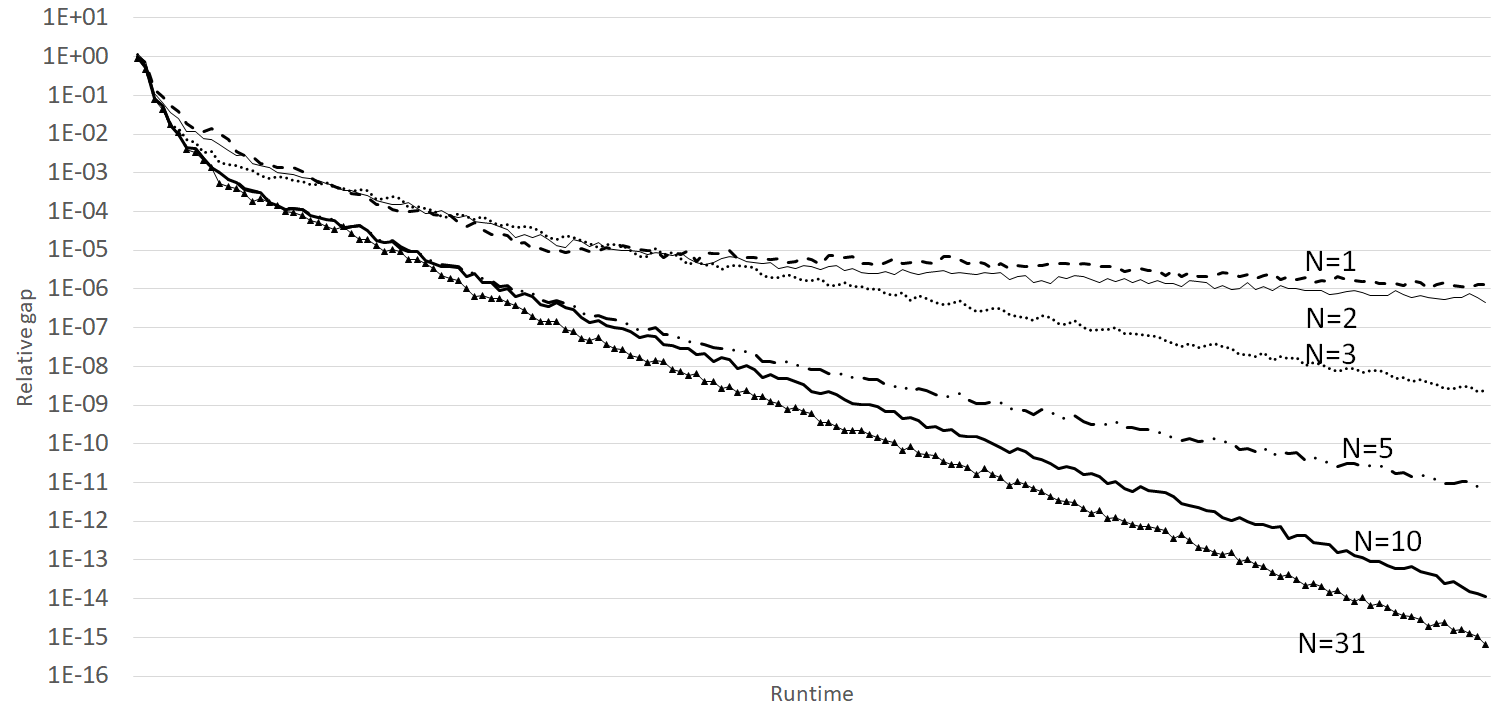}}
\subfigure[Convergence for A-TAP GP with $N$-link cost dependency]{\label{fig:CS_b}\includegraphics[width=0.8\textwidth]{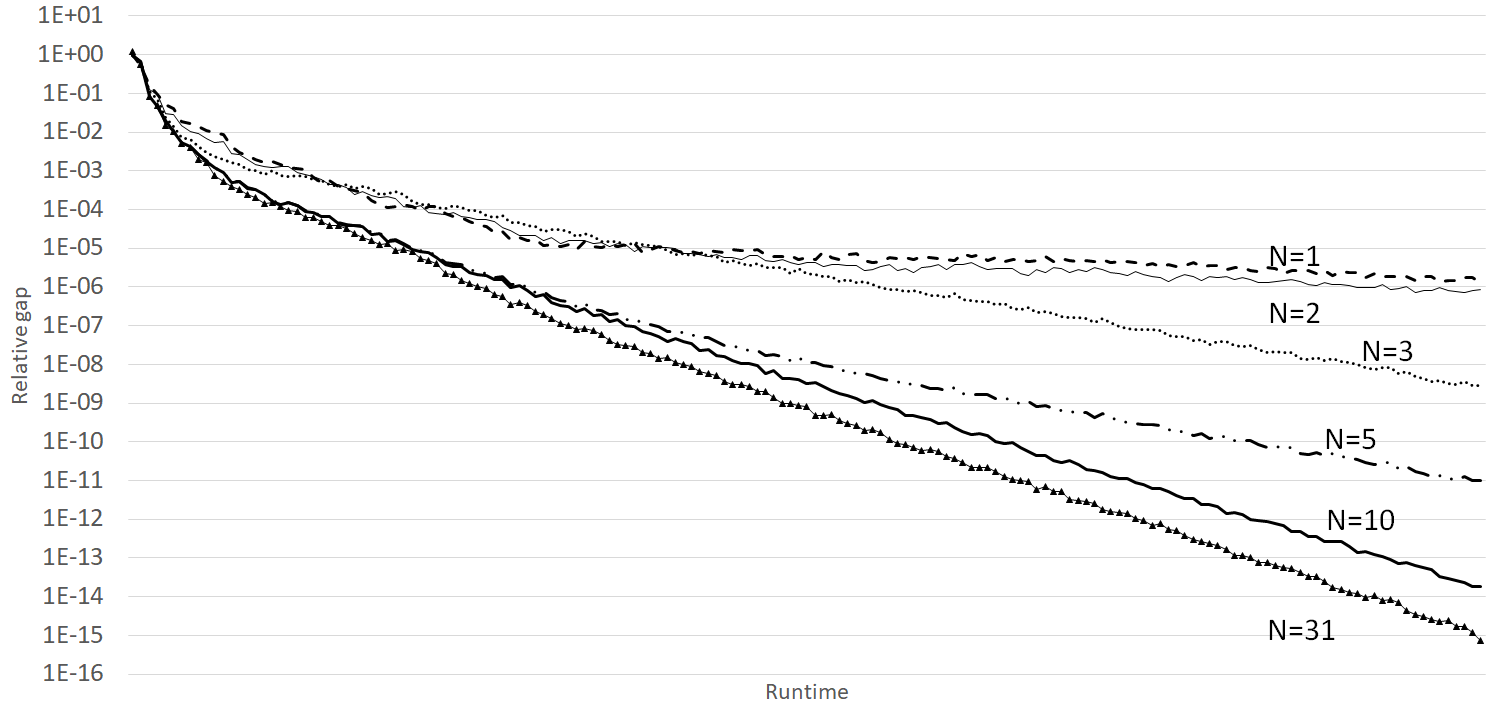}}
\subfigure[Asymmetric to symmetric weight matrix convergence]{\label{fig:CS_c}\includegraphics[width=0.8\textwidth]{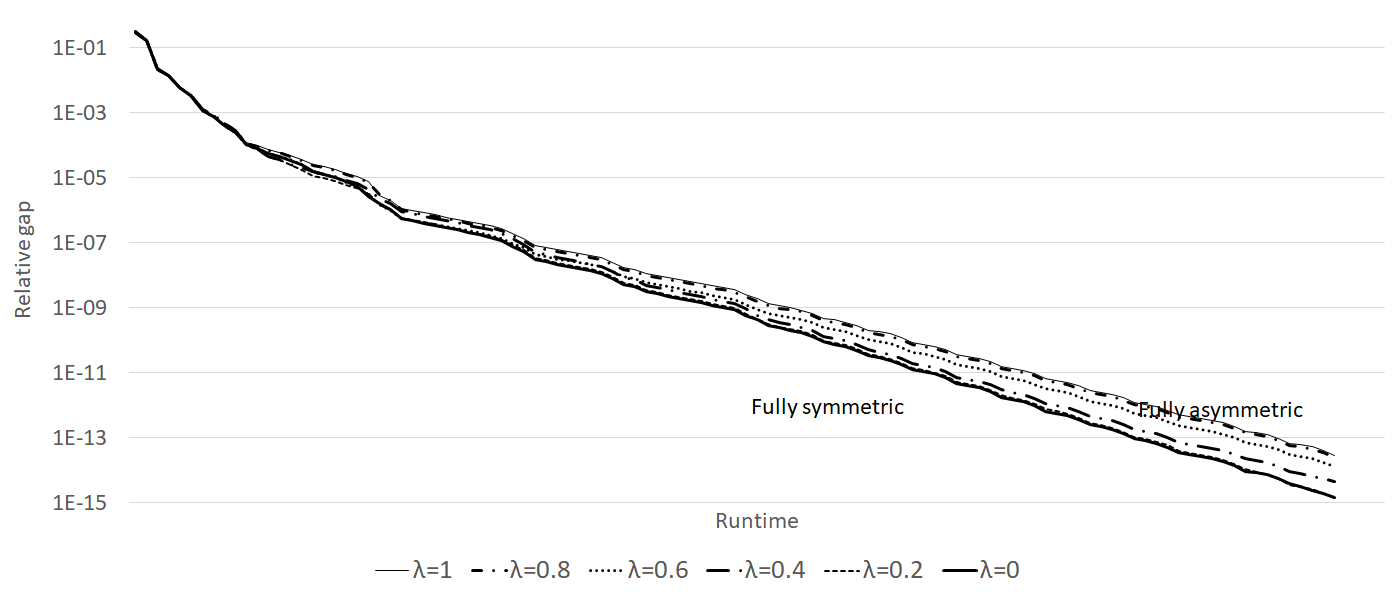}}
\caption{Experimental results for Chicago-sketch network}
\label{fig:CS_compiled_results}
\end{figure}

\clearpage

\begin{figure}
\centering
\subfigure[Convergence for S-TAP GP with $N$-link cost dependency]{\label{fig:Barca_a}\includegraphics[width=0.8\textwidth]{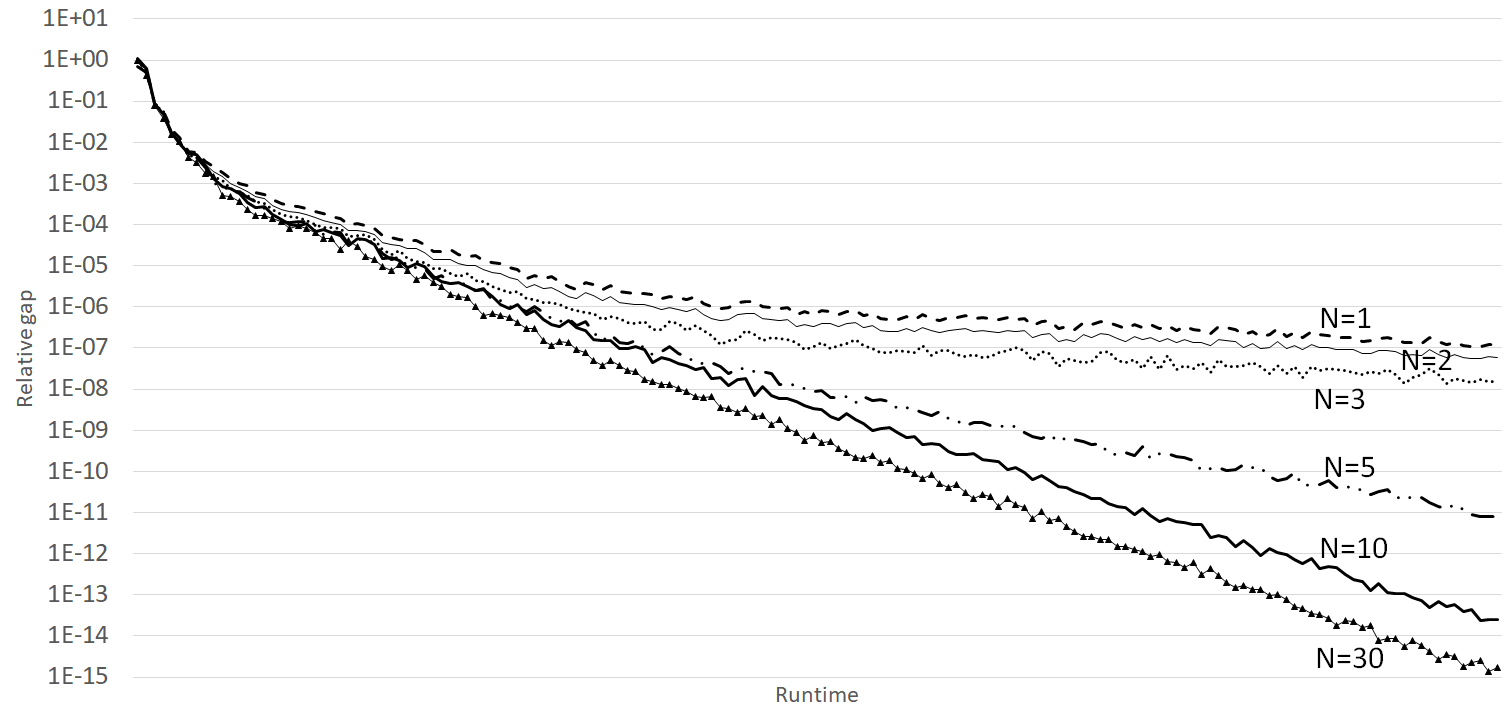}}
\subfigure[Convergence for A-TAP GP with $N$-link cost dependency]{\label{fig:Barca_b}\includegraphics[width=0.8\textwidth]{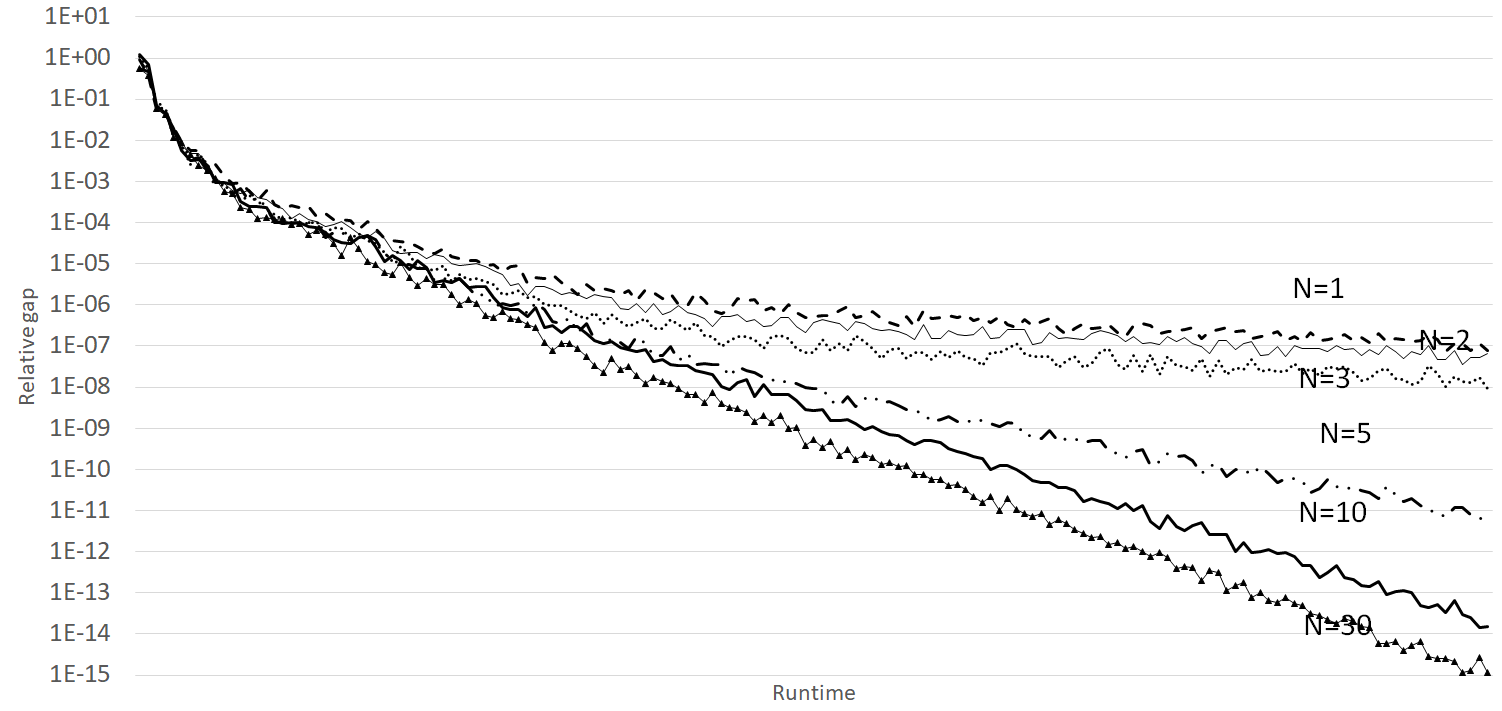}}
\subfigure[Asymmetric to symmetric weight matrix convergence]{\label{fig:Barca_c}\includegraphics[width=0.8\textwidth]{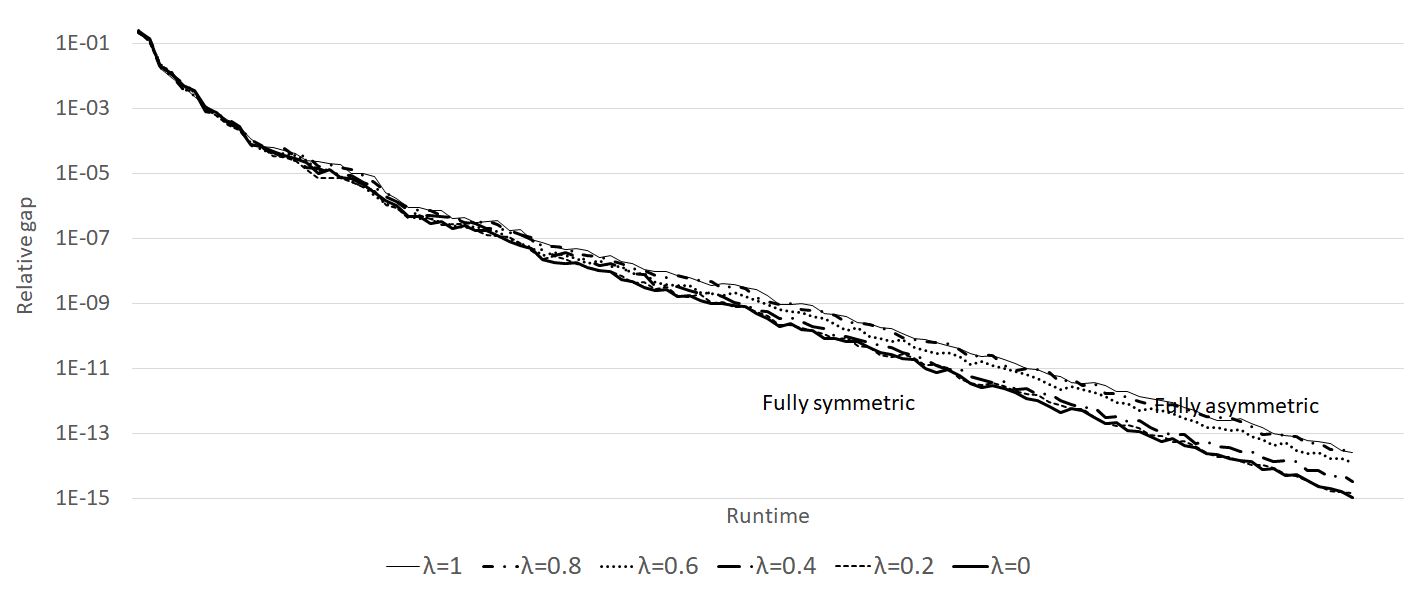}}
\caption{Experimental results for Barcelona network}
\label{fig:Barca_compiled_results}
\end{figure}

\clearpage

\begin{figure}
\centering
\subfigure[Convergence for S-TAP GP with $N$-link cost dependency]{\label{fig:CR_a}\includegraphics[width=0.8\textwidth]{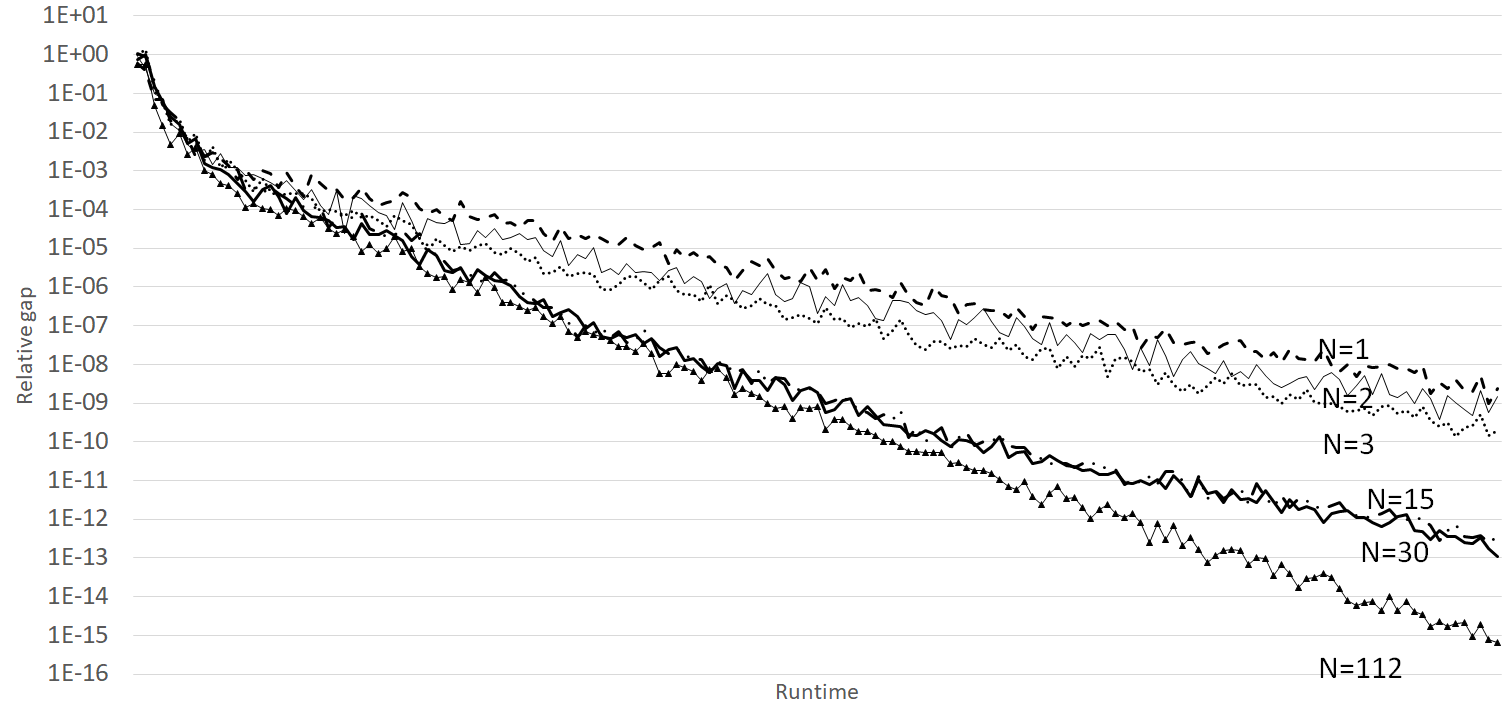}}
\subfigure[Convergence for A-TAP GP with $N$-link cost dependency]{\label{fig:CR_b}\includegraphics[width=0.8\textwidth]{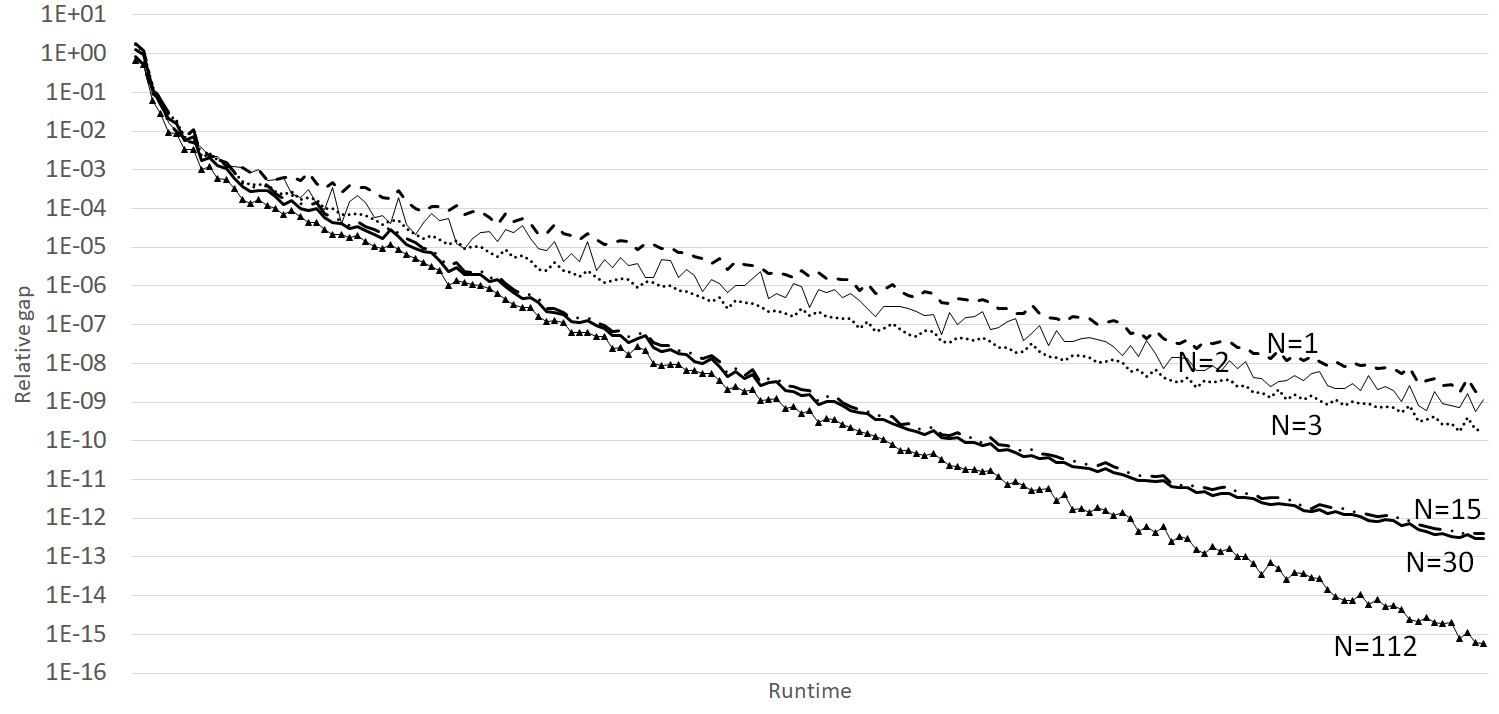}}
\subfigure[Asymmetric to symmetric weight matrix convergence]{\label{fig:CR_c}\includegraphics[width=0.8\textwidth]{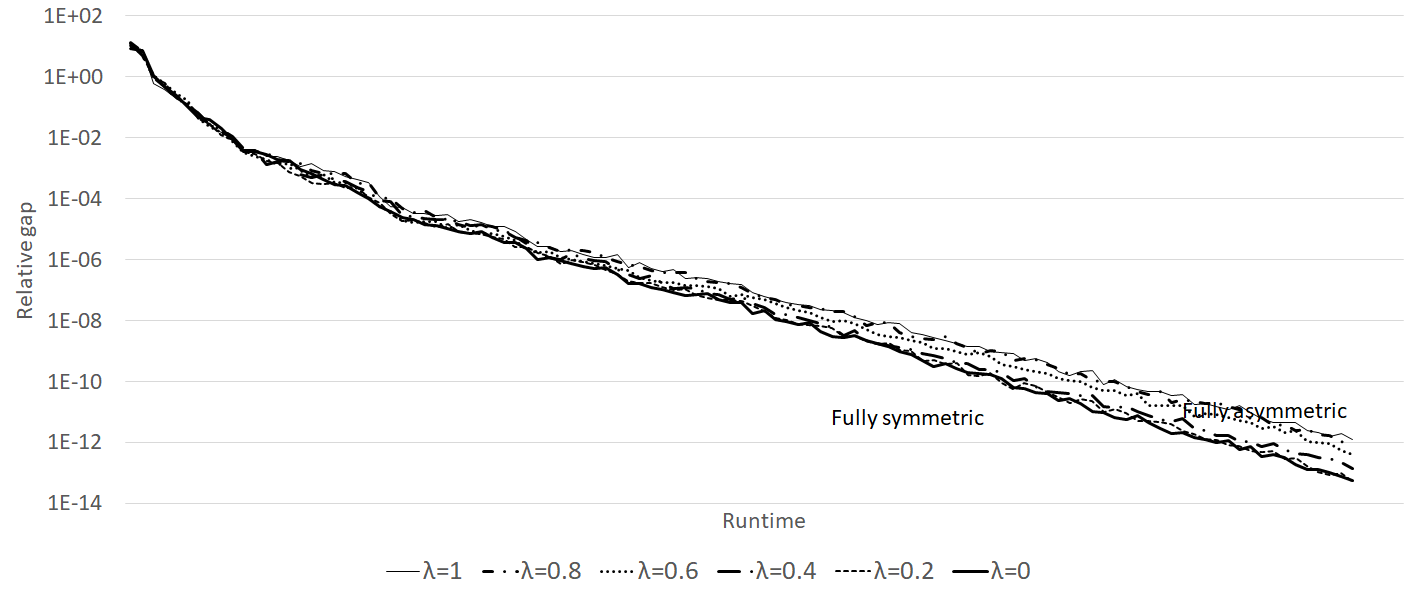}}
\caption{Experimental results for Chicago-regional network}
\label{fig:CR_compiled_results}
\end{figure}

\clearpage

\begin{figure}
\noindent\makebox[\textwidth]{%
\includegraphics[scale=0.85]{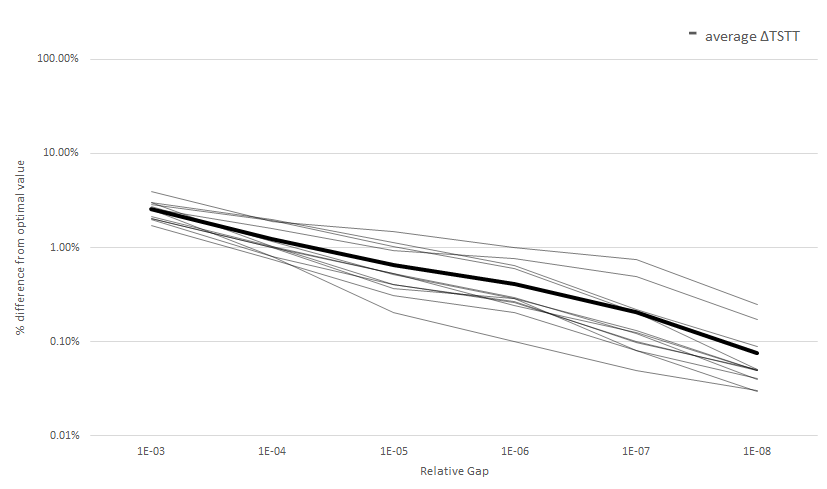}}
\centering
\caption{$\Delta$TSTT trends for different gap levels}
\label{fig:TSTT_trend}
\centering
\end{figure}

\clearpage

\begin{figure}
\noindent\makebox[\textwidth]{%
\includegraphics[scale=0.85]{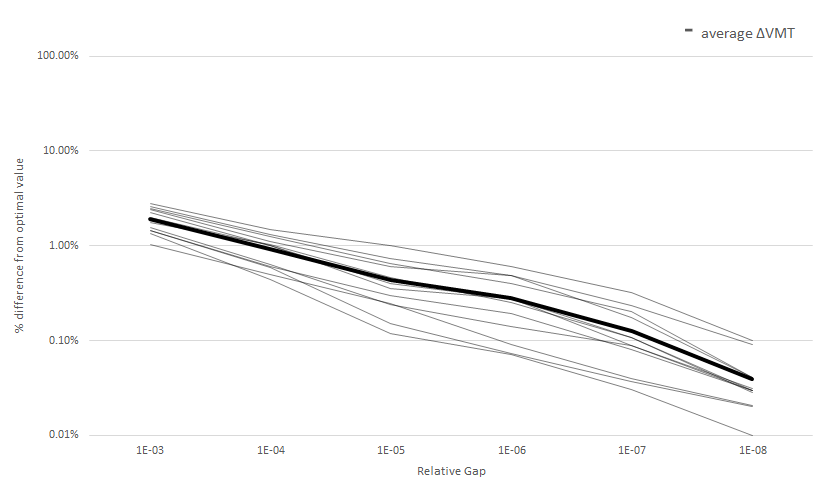}}
\centering
\caption{$\Delta$VMT trends for different gap levels}
\label{fig:VMT_trend}
\centering
\end{figure}

\clearpage

\begin{figure}
\noindent\makebox[\textwidth]{%
\includegraphics[scale=0.85]{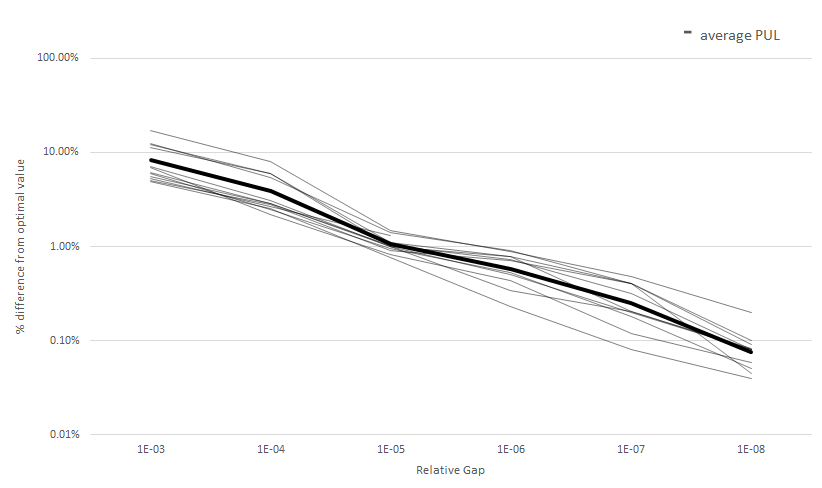}}
\centering
\caption{PUL trends for different gap levels}
\label{fig:PUL_trend}
\centering
\end{figure}

\clearpage

List of figure captions:
\begin{itemize}
\item Figure 1. S-TAP example with multiple extreme points
\item Figure 2. Toy Network
\item Figure 3. Condition number behavior of Sioux Falls problem instances
\item Figure 4. Experimental results for Sioux Falls network
\item Figure 5. Sioux Falls asymmetric to symmetric weight matrix convergence
\item Figure 6. Experimental results for Eastern Massachusetts network
\item Figure 7. Experimental results for Chicago-sketch network
\item Figure 8. Experimental results for Barcelona network
\item Figure 9. Experimental results for Chicago-regional network
\item Figure 10. $\Delta$TSTT trends for different gap levels
\item Figure 11. $\Delta$VMT trends for different gap levels
\item Figure 12. PUL trends for different gap levels

\end{itemize}

\end{document}